\documentclass[11pt,a4paper]{amsart}
\usepackage{amsmath,mathtools,mathrsfs,multirow}

\usepackage[scr=boondox]{mathalpha}

\calclayout
\usepackage{hyperref}
\usepackage{stmaryrd}
\usepackage{enumerate}
\usepackage{comment}
\usepackage{graphicx}
\usepackage{amsfonts}
\usepackage{amssymb}
\usepackage{wrapfig,lipsum,booktabs}
\usepackage{placeins}
\usepackage[hang]{subfigure}
\usepackage[utf8]{inputenc}
\usepackage[english]{babel}
\usepackage{amsthm}
\usepackage{subfigure} 
\usepackage[utf8]{inputenc}
\usepackage[margin=0.95in]{geometry} 
\usepackage[final]{pdfpages}
\newcounter{count}\numberwithin{count}{section}
\newtheorem{theorem}[count]{Theorem}
\newtheorem*{theorem*}{Theorem}
\newtheorem*{remark}{Remark}
\newtheorem{lemma}[count]{Lemma}
\theoremstyle{definition}

\theoremstyle{notation}

\theoremstyle{corollary}
\newtheorem{corollary}[count]{Corollary}

\newtheorem{Prop}[count]{Proposition}
\usepackage{times}
\usepackage{blindtext}
\usepackage{setspace}
\author{Yazan Alamoudi}
\address{Department of Mathematics\\
University of Florida\\
Gainesville\\
FL 32611\\
United States
}
\email{yazanalamoudi@ufl.edu}

\dedicatory{Dedicated to my mentor Krishnaswami Alladi, in appreciation of his guidance and insight}

\title{On subradically sifted sums related to Alladi's higher order duality between prime factors}

\begin{document}

\begin{abstract}
In this paper, I utilize a variant of the Selberg--Delange method to find quantitative estimates of the sums \[M_{k,\omega}(x,y)=\sum_{\substack{p_{1}(n)> y\\ n\leq x} } \mu(n) {\omega(n)-1\choose k-1},\] where $y$ can grow with $x$ but we must have $y\leq Y_0\exp(\mathscr{p}\frac{\log x}{(\log\log (x+1))^{1+\epsilon}})$ with $Y_0,\mathscr{p},\epsilon>0$. Moreover, I give preliminary upper bounds for the general range $1.9\leq y\leq x^{\frac{1}{k}}$. In addition, I formalize the notions of subradical and radical dominance and discuss their relevance to the analytic approach to the study of arithmetic functions. Lastly, I give a fascinating formula related to the derivatives of the gamma function and the Hankel contour, which should be relevant for those employing the Selberg--Delange method to obtain higher-order terms.
\end{abstract}

\maketitle

\section{Introduction}
More than four and a half decades ago, K. Alladi \cite{KA77} stated his four famous duality relations. Most relevant for this paper are the last two of these relations. Specifically, let\footnote{For the remainder of this paper, we will abide by the convention $p_1(1)=\infty$ in line with \cite{KA82}.} $p_k(n),P_k(n)$ denote the $k$-th smallest and $k$-th largest prime factors, respectively. Alladi's third and fourth duality relations, found in \cite{KA77}, state that if $f(1)=0$ then
\[\sum_{d|n} \mu(d) {\omega(d)-1\choose k-1}f(p_1(d))=(-1)^kf(P_k(n)),\]
and
\[\sum_{d|n} \mu(d) {\omega(d)-1\choose k-1}f(P_1(d))=(-1)^kf(p_k(n)).\]
The case $k=1$ was used to find the intriguing infinite sum 
\[\sum_{p_1(n)\equiv l \pmod{m}}\frac{\mu(n)}{n}=\frac{-1}{\phi(m)}\]
along with other quantitative and arithmetic results. Alladi's first-order duality has gotten much attention over the years. Alladi himself revisited the topic again in \cite{KA82} by studying the sum 
\[M(x,y)=\sum_{\substack{p_{1}(n)> y\\ n\leq x} } \mu(n). \]
In view of the relation,
\[\sum_{n\leq x}\mu(n)f(p_1(n))=-\sum_{p\leq x}f(p)M(\frac{x}{p},p),\]
Alladi described (see \cite{KA82}) $M(x,y)$ as the ``basic computational" tool for studying the sums  
\[M_f(x)=\sum_{n\leq x}\mu(n)f(p_1(n)).\] On the other hand, the importance of estimating the size of the above sums is demonstrated again and again in Alladi's original paper \cite{KA77}. 

More recently, Alladi and Johnson \cite{KAJJ24} expanded the horizon by studying the case $k=2$ and demonstrating several results, including
\[\sum_{p_1(n)\equiv l \pmod{m}}\frac{(\omega(n)-1)\mu(n)}{n}=\frac{1}{\phi(m)}.\]
This movement maintains momentum, as Alladi and Sengupta (see \cite{KASS,KAJJ24}) are handling the general relation
\[\sum_{p_1(n)\equiv l \pmod{m}}\frac{ {\omega(n)-1\choose k-1}\mu(n)}{n}=\frac{(-1)^k}{\phi(m)}.\]
Although the basic study of the $k$-th order duality, for $k$ greater than one, has already commenced, the sums 
\[M_{k,\omega}(x,y)=\sum_{\substack{p_{1}(n)> y\\ n\leq x} } \mu(n) {\omega(n)-1\choose k-1},\tag{1.1}\]
have not been studied in the aforementioned projects.\footnote{However, K. Alladi and I are currently preparing \cite{YAKASS} a comprehensive treatment for the case $k=2$.}\\

The primary purpose of this paper is to give refined quantitative estimates for (1.1) with general $k$ in the cases $y\leq Y_0\exp(\mathscr{p}\frac{\log x}{(\log\log (x+1))^{1+\epsilon}})$ with $Y_0,\mathscr{p},\epsilon>0$. In what follows, I state the main theorem which features a term defined by $\Gamma_{N+1,m}:=\lim_{z\to N+1}(\frac{d}{dz})^m\frac{1}{\Gamma(-z)}<\infty$, which will be shown to be convergent and given an explicit formula in the next section.

\begin{theorem}
Let $N,k-1\in \mathbb{N}$ and $Y_0,\mathscr{p},\epsilon\in\mathbb{R}_{>0}$. If $y>x^{\frac{1}{k}}$, then $M_{k,\omega}(x,y)=(-1)^{k-1}$. Otherwise, there is a constant $C_{N,k-1}$, that may depend on $N$, $k-1$, $Y_0$, $\epsilon$ and $\mathscr{p}$ but is otherwise absolute (in particular, $C_{N,k-1}$ does not depend on $x$ or $y$), such that whenever $1.9\leq y\leq\min{\{ Y_0\exp(\mathscr{p}\frac{\log x}{(\log\log (x+1))^{1+\epsilon}}),x^{\frac{1}{k}}\}}$, the following holds. \[M_{k,\omega}(x,y)=\frac{x}{\log x}\sum_{\substack{1\leq i'\leq i\leq N \\1\leq j\leq J\leq k-1\\0\leq j'\leq k-1} }\mu_{N,k,i',i,j',j,J}(\log\log y)^{j'}(\log y)^{i'}\frac{\Gamma_{i,j}}{(\log x)^i}(\log{\log x})^{J-j}+\varepsilon_{N,k-1}(x,y) \tag{1.2},\]
and 
\[M_{k,\omega}(x,y)=\frac{x}{\log x}\sum_{\substack{1\leq i\leq N \\1\leq J\leq k-1\\1\leq j\leq J} }\phi_{i,J}(y)\frac{\Gamma_{i,j}}{(\log x)^i}(\log{\log x})^{J-j}{J\choose j}  + \varepsilon_{N,k-1}'''(x,y),\tag{1.3}\]
with
\[|\varepsilon_{N,k-1}(x,y)|\leq C_{N,k-1}\frac{x(\log\log (x+1))^{k-1}}{\log x}\bigg\{(\frac{\log y}{\log x})^{N+1}+\frac{1}{\log x\exp (c''\sqrt{\log y})}\bigg\},\]
and 
\[|\varepsilon'''_{N,k-1}(x,y)|\leq C_{N,k-1}\frac{x(\log\log (x+1))^{k-1}}{\log x}(\frac{\log y}{\log x})^{N+1},\]
where $\phi_{i,J}(y)$ are some functions of $y$ given an explicit formula, in terms of some constants, later in the manuscript\footnote{See equations (4.8), (4.9), and (4.16).} such that 
\[\phi_{i,J}(y)\ll(\log y)^{i+1}(\log \log y)^{k-1}.\]
    
\end{theorem}
\begin{remark}
K. Alladi has pointed out to me that if $y$ is fixed then (1.2) only gives an upper bound. In view of his remark, I have decided to include\footnote{Equation (1.3) was initially only used to prove (1.2).} (1.3) in the main theorem. Note that (1.3) actually gives for arbitrary but fixed $y_0$
\[M_{k,\omega} (x,y_0)=\frac{x}{\log x}\sum_{\substack{1\leq i\leq N \\1\leq J\leq k-1\\1\leq j\leq J} }\phi_{i,J}(y_0)\frac{\Gamma_{i,j}}{(\log x)^i}(\log{\log x})^{J-j}{J\choose j}  + O_{N,k-1,y_0}(\frac{x(\log\log x)^{k-1}}{(\log x)^{N+2}}).\]
On that same note, it is worth mentioning that K. Alladi and S. Sengupta have recently\footnote{It is important to emphasize that the results of this paper were obtained independently by the author. In particular, at the time of the initial submission and first posting on arXiv, \cite{KASS} was still in preparation and I did not have access to the full manuscript; I had only seen the treatment of $M_\omega(x)$, a closely related variant of $M_{2,\omega} (x,1)$, which was consistent with results I had obtained prior to seeing it.} obtained sharp estimates for $M_{\omega^k}(x)=\sum_{n\leq x} \mu(n)\omega(n)^k$ in the form of a finite sum featuring decreasing powers of $\log x$ (see \cite{KASS}, where the method for $k=2$ is completely elementary and the general $k$ is stated). 
\end{remark}

Observe that the order of the upper bound $O(\exp(\mathscr{p}\frac{\log x}{(\log\log (x+1))^{1+\epsilon}}))$ for the sifting is considerably larger than the typical $O(\exp{c\sqrt{\log x}})$ or $O(\exp{c\sqrt[n]{\log x}})$ for some $n\in\mathbb{N}\setminus\{0,1\}$. Actually, one can argue that there is a more categorical difference between the bounds of the main theorem and the more typical sifting ranges found by the Selberg--Delange method. For this, allow me to introduce some terminology.

I will say $f$ is subradically dominated by $g$ and write $f<_{\sqrt[\forall]{\,\,}}g$ if for any $\epsilon > 0$ $f(x)=o(g(x)^{\epsilon})$. Functions subradically dominated by $x$ are simply called subradical. On the other hand, I will say that $f$ is radically dominated by $g$, and write $f<_{\sqrt[\exists]{\,\,}}g$, if there exists some $0<\epsilon < 1$ such that $f=o(g^{1-\epsilon})$. The key thing to point out is that although both the sifting bound here and the typical ones are subradical, the bound here subradically dominates the more classical bounds.\footnote{In particular, for any $\epsilon',\epsilon>0$ we have \(e^{(\log{x})^{1-\epsilon'}} <_{\sqrt[\forall]{\,\,}}\exp(\mathscr{p}\frac{\log x}{(\log\log (x+1))^{1+\epsilon}})<_{\sqrt[\forall]{\,\,}}x.\)}

Unfortunately, the set of subradical functions has neither a maximal element for the order $<_{\sqrt[\forall]{\,\,}}$ nor one for the order $<_{\sqrt[\exists]{\,\,}}$.\footnote{In Section 6, I prove that there is no such maximal element and that the set of subradical functions actually forms a ring.} In other words, if you do not get all subradical sieves, there will always be sieves that subradically dominate your sieve.\footnote{This is not necessarily to be expected as the set of subradical functions has multiple least upper bounds in the set of all functions. Specifically, for $<_{\sqrt[\forall]{\,\,}}$, all of $x$, $\frac{x}{\log x}$, and $\sqrt{x}$ are, although not subradical, least upper bounds for the set of subradical functions.} On the other hand, if you look at the log of the sifting variable's upper bound, then it is in fact true that for any $\epsilon,\epsilon'>0$ we have that ${\log{\exp{\frac{\log x}{(\log\log x)^{1+\epsilon}}}}=\frac{\log x}{(\log\log x)^{1+\epsilon}}}$ is maximal for all subradical functions with respect to the strict partial order $<_{\sqrt[\exists]{\,\,}}$ whereas ${\log e^{(\log{x})^{1-\epsilon'}}\!\!={(\log{x})^{1-\epsilon'}}}$ is not. From this perspective, it can be argued that if one can't get all subradical functions, the sifting bounded by $Y_0\exp(\mathscr{p}\frac{\log x}{(\log\log (x+1))^{1+\epsilon}})$ is an optimal alternative.

The above argument may not convince every reader. After all, this paper introduced the strict partial orders $<_{\sqrt[\exists]{\,\,}}$ and $<_{\sqrt[\forall]{\,\,}}$, and used them to justify the sifting rates. This is a valid concern, but I still invite the reader to reflect on the two orders, as they encode a common pattern in the typical error terms that arise in analytic number theory. Nonetheless, let us set aside the orderings and provide a more concrete justification for the significance of the improved sifting rate. As early as $k=3$, if we let $\epsilon=\frac{1}{10^{10}}$, then even with $Y_0=1$ and $\mathscr{p}=2$, we have $\sqrt[3]{x}<\exp(2\frac{\log x}{(\log\log (x+1))^{1+\epsilon}})$ for $x<e^{e^{5.9}}-1$ which exceeds the accepted estimate of the number of protons of the observable universe (see \cite{UNIPTN}). Given Alladi's characterization of $M(x,y)$ as a ``basic computational tool" in \cite{KA82}, I see this as meeting a convincing benchmark for such purposes.

To conclude the introduction, I will briefly outline the manuscript's structure. The next section will give a basic setup and go over preliminaries. The third handles errors related to the sifting variable. The fourth will simply put everything together, as the most technical work is largely finished, and complete the proof of all the main quantitative statements. The fifth will give preliminary upper bounds for sifting rates beyond $\exp(\mathscr{p}\frac{\log x}{(\log\log (x+1))^{1+\epsilon}})$. The paper concludes its technical matter with section 6, which exposits a general framework for the subradical and radical ordering of functions. Although only specific instances of these notions are directly related to the main theorem, I believe that these notions encode a materially significant boundary in certain topics in analytic number theory, and more attention may be a valuable investment.

\section{An overview of the approach}
As is typical with the Selberg--Delange method (see \cite{KAQJ}), we will first consider\footnote{Parts of this paper grew out of material in my oral exam and my contributions to \cite{YAKASS}.}
\[
   f(s,y,z)=\sum_{\substack{n=1\\p_{1}(n)> y} }^\infty \frac{\mu(n)^2 z^{\omega(n)-1}}{n^s}=z^{-1}\prod_{p> y}\bigl(1+\frac{z}{p^s}\bigr)=g(s,y,z)\zeta(s)^z, \tag{2.1}
\]

\noindent where 

\[
  g(s,y,z)=z^{-1}\left(\prod_{p}\bigl(1+\frac{z}{p^s}\bigr)\bigl(1-\frac{1}{p^s}\bigr)^z\right)\prod_{p\leq y}\bigl(1+\frac{z}{p^s}\bigr)^{-1}\tag{2.2}.
\]

To begin with, note that 

\[\frac{(-1)^k}{k-1!}(D_z)^{k-1}f(s,y,z)\mid_{z=-1}=\sum_{\substack{n=1\\p_{1}(n)> y} }^\infty \frac{\mu(n)^2 {\omega(n)-1\choose k-1}(-1)^{\omega(n)}}{n^s}=\sum_{\substack{n=1\\p_{1}(n)> y} }^\infty \frac{\mu(n) {\omega(n)-1\choose k-1}}{n^s}.\]

With the above in mind, we consider a basic consequence of the well-known truncated Perron formula.
\begin{lemma}
For $x,y\geq 1.9$ and $k\geq1$ we have 
\[
   \bigg|\frac{1}{2\pi i}\int_{a-iT}^{a+iT}\frac{(-1)^k}{k-1!}(D_z)^{k-1}f(s,y,z)\mid_{z=-1}\frac{x^s}{s}ds - \sum_{\substack{p_{1}(n)> y\\ n\leq x} } \mu(n) {\omega(n)-1\choose k-1}\bigg|\leq e\times\frac{x(5\log x)^{2k}}{T},\tag{2.3}
\]
and 
\[
  \hspace*{-1cm} \bigg|\frac{1}{2\pi i}\int_{a-iT}^{a+iT}\frac{(-1)^k}{k-1!}(D_z)^{k-1}f(s,y,z)\mid_{z=-1}\frac{x^s}{s}ds - \sum_{\substack{p_{1}(n)> y\\ n\leq x} } \mu(n) {\omega(n)-1\choose k-1}\bigg|\leq  e\times\frac{x}{T}\times\frac{1}{k-1!}(D_z)^{k-1}f(a,y,z)\mid_{z=1}.\tag{2.4}
\]
   
Where $a=1+\frac{1}{\log{x}}$ and $T$ is some suitable function of $x$ to be specified.

 \end{lemma}

\begin{proof}

Both (2.3) and (2.4) follow from instances of the truncated Perron formula with the well-known explicit error (see, for example, \cite{Hild}) as 
\[R(T)\leq \frac{x^a}{T}\sum_{\substack{n=1\\p_{1}(n)> y} }^\infty \frac{|\mu(n) {\omega(n)-1\choose k-1} |}{n^a}\leq e\times\frac{x}{T}\times\frac{1}{k-1!}(D_z)^{k-1}f(a,y,z)\mid_{z=1}
\]\[\leq\frac{x^a}{T}\sum_{n=1}^\infty \frac{\omega(n)^{k-1}}{n^a}\leq\frac{x^a}{T}\bigg(\zeta(a)(\sum_p\frac{1}{p^a})\bigg)^k\leq e \times \zeta(a)^{2k} \times\frac{x}{T}.\]
Note that, for real $\sigma>0$, we have 
\[\frac{1}{\sigma-1}\leq\zeta(\sigma)\leq \frac{\sigma}{\sigma-1}.\]
It follows that
$$ e \times \zeta(a)^{2k} \times\frac{x}{T}\leq e\times\frac{x(5\log x)^{2k}}{T}$$
since $1+\frac{1}{\log x}\leq 1+\frac{1}{\log 1.9}<5$. 
\end{proof}
Now, we proceed by noting\footnote{A heuristic use of differentiating with respect to $z$ appears in \cite{KAJJ24}. The main argument and results in this paper are rigorous and self-contained (given standard background and the cited results).} that the Leibniz rule gives
\[\frac{(-1)^k}{k-1!}(D_z)^{k-1}f(s,y,z)\mid_{z=-1}=\frac{(-1)^k}{k-1!\zeta(s)}\sum_{i=0}^{k-1}{k-1\choose i}(\log\zeta(s))^{k-i-1}(D_z)^{i}g(s,y,z)\mid_{z=-1}.\tag{2.5}\]

With this in mind, we will obtain the main theorem by adopting the classic (see, for example, \cite{MV}) contour shift with a slight variant of a larger choice of $T$. Specifically, in Section 4, I will replace the contour in (2.3) with
\begin{enumerate}[(i)]
    \item horizontal lines connecting $a\pm iT$ to $1-\frac{c}{\log T}\pm iT$,
    \item vertical lines connecting $1-\frac{c}{\log T}\pm iT$ to $1-\frac{c}{\log T}\pm \frac{i}{\log{x}}$, and
    \item the important contour $\mathcal{C}(x)$ composed of horizontal lines connecting  $1-\frac{c}{\log T}\pm \frac{i}{\log{x}}$ to $1\pm \frac{i}{\log{x}}$ and the semi circle $\{1+ \frac{e^{i\theta}}{\log{x}}:-\frac{\pi}{2}\leq\theta\leq\frac{\pi}{2}\}$.
\end{enumerate}
Controlling $\mathcal{G}_i(s,y,-1)=(D_z)^{i}g(s,y,z)\mid_{z=-1}$ as well as its $s$-derivatives is among the most technically involved aspects and is the content of the next chapter. On the other hand, an invocation of the Hankel formula is often standard in such investigations. Since the approach here involves taking derivatives and truncating past the first term, the reader will be correct to expect that we will need a more general formula related to the derivatives of $\Gamma(z)$.

This is precisely the content of the upcoming lemma.\footnote{After the first submission, I learned of a related formula in \cite{TR}. The results here were obtained independently, by a different method.}

\begin{lemma}
Let $\mathcal{H}(x)$ be the transformed region corresponding to the region $\mathcal{C}(x)$ under the standard substitution $s=1+\frac{w}{\log x}$ and let $\mathcal{H}=\lim_{x\to \infty}\mathcal{H}(x)$ be the resulting Hankel contour. Set $\zeta(1):=\gamma$. Then

$$\frac{1}{2\pi i}\int_{\mathcal{H}}\log(w)^mw^{N+1} e^w dw=\lim_{z\to N+1}(\frac{d}{dz})^m\frac{1}{\Gamma(-z)}:=\Gamma_{N+1,m}$$
$$=(-1)^N\sum_{\text{$i$ is odd } i+j+k=m}(-1)^\frac{i-1}{2}k! {m\choose i,j,k}\pi^{i-1}Y_{j}(-0!\zeta(1),1!\zeta(2),\dots,(j-1)!(-1)^j\zeta(j)) c(N+2,k+1)$$

and
\[\frac{1}{2\pi i}\int_{\mathcal{H}(x)}\log(w)^mw^{N+1} e^w dw=\Gamma_{N+1,m}+O_{N+1,m}(\exp{\frac{-c^{-}\log x}{\log T}}),\tag{2.6}\]
with $c^{-}$ is any positive real number with $c^{-}<c$, $Y_j$ the exponential Bell polynomial, and $c(m,n)$ the unsigned Stirling number of the first kind. Note that, in particular $\Gamma_{N,m}\in\mathbb{Q}(\gamma,\zeta(2),\dots,\zeta(m-1))$.
\end{lemma}
 
\begin{proof}
The first equality follows directly from differentiation under the integral, which is justified by the dominated convergence theorem. On the other hand, the equality of $\Gamma_{N+1,m}$ to the sum in the lemma can be established independent of the Hankel contour by relating the quantity to well-known combinatorial entities. More specifically, begin by noting the well-known fact in elementary combinatorics stating\footnote{A standard way to prove (2.7) is to interpret LHS as a value of an elementary symmetric function and compare with the recurrences of $c(m,n)$. Alternatively, one may use a combination of Taylor's theorem and the relevant generating polynomial.} \[\frac{1}{k!}(\frac{d}{dz})^kz(z-1)\dots (z-N)\mid_{z=N+1}= c(N+2,k+1).\tag{2.7}\] To complete the argument, note that, in view of the results for $\Gamma^{(n)}(1)$ stated in \cite{Donal} in terms of the exponential Bell polynomials and the well-known relation \[\frac{1}{\Gamma(-z)}=-\frac{\sin(\pi z)\Gamma(1+z)}{\pi}=-\frac{\sin(\pi z)\Gamma(z-N)z(z-1)\dots (z-N)}{\pi},\] an application of the Leibniz rule gives
$$\lim_{z\to N+1}(\frac{d}{dz})^m\frac{1}{\Gamma(-z)}=-\lim_{z\to N+1}\sum_{i=0}^m {m\choose i}(\frac{d}{dz})^{i}\frac{\sin(\pi z)}{\pi}(\frac{d}{dz})^{m-i}\Gamma(z-N)z(z-1)\dots (z-N)$$
$$=(-1)^N\lim_{z\to N+1}\sum_{\text{$i$ is odd } 0\leq i\leq m}(-1)^\frac{i-1}{2} {m\choose i}\pi^{i-1}(\frac{d}{dz})^{m-i}\Gamma(z-N)z(z-1)\dots (z-N)$$
$$=(-1)^N\sum_{\text{$i$ is odd } i+j+k=m}(-1)^\frac{i-1}{2}k! {m\choose i,j,k}\pi^{i-1}\Gamma^{(j)}(1)c(N+2,k+1)$$
$$=(-1)^N\sum_{\text{$i$ is odd } i+j+k=m}(-1)^\frac{i-1}{2}k! {m\choose i,j,k}\pi^{i-1}Y_{j}(-0!\zeta(1),1!\zeta(2),\dots,(j-1)!(-1)^j\zeta(j)) c(N+2,k+1).$$
On the other hand, once the first formula is established, proving (2.6) is completely standard. First, note that $\mathcal{H}=\mathcal{H}(x)\cup\mathcal{H}_\pm$ where $\mathcal{H}_\pm=\{w=u\pm i :-\infty\leq u \leq \frac{c\log x}{\log T}\}.$ From here we conclude that
\[\hspace*{-0.7cm}\int_{\mathcal{H}_\pm}\log(w)^mw^{N+1} e^w dw\ll_{N,m}\int_{\frac{c\log x}{\log T}}^\infty u^{N+2} e^{-u} dw\ll_{N,m} \biggl|\frac{c\log x}{\log T}\biggr|^{N+3}\exp{\frac{-c\log x}{\log T}}\ll_{N,m}\exp{\frac{-c^{-}\log x}{\log T}},\]
for any $c^{-}<c$.
\end{proof}
\begin{remark}What I find fascinating in the above is that $\Gamma_{N,m}\in\mathbb{Q}(\gamma,\zeta(2),\dots,\zeta(m-1))$ and these constants appear in Theorem 1.1 which relate to $k$-th prime factors. From this perspective, it is natural to think that if the common current attitude\footnote{The exact transcendence status of many of these constants remains unknown to the best of my knowledge.} in the mathematical community is true, then these constants have progressively larger transcendence degrees. In this context, I wonder, from a naive place, why asking questions about higher-order prime factors would come at the cost of ``infinite" complexity?\end{remark}
At this point, we switch to the much more technical task of bounding $(D_s)^j(D_z)^{j'}g(s,y,z)\mid_{z=-1}$, by an elementary function\footnote{One can extend the range of $y$ by bounding $\mathcal{G}_i(s,y,-1)$ by an elementary function of $y$ and $T$. However, this requires a bit more work and is unnecessary for the main result (Theorem 1.1). Hence, it was avoided in the current manuscript.} of $y$. The reader should note that the purpose of Section 3 is to establish Theorem 3.1, as this is the only result that will be used in the proofs in Sections 4 and 5. As such, the technical details arising in the next section will not be revisited later in the paper. However, Section 4 ends with an interesting remark building on the remarks from Section 3; this remark is not used elsewhere. 
\section{Controlling the error related to the sifting variable}  
The main purpose of this section is to establish:
\begin{theorem}

There exist constants $m_{N,k-1}$ (that do not depend on $y$ or $T$) such that for $\sigma\geq1-\frac{c}{\log T}$, $1.9\leq y\leq Y_0T^{\mathscr{p}}$ with $T\geq T_0$, $j\leq N$, and $j'\leq k-1$ we have
\[\left|(D_s)^j(D_z)^{j'}g(s,y,z)\mid_{z=-1}\right|<m_{N,k-1}(\log y)^{j+1}(\log\log (y+1))^{j'}.\tag{3.1}\] 
In addition, for $N,k-1\geq0$ we have \[(D_s)^N(D_z)^{k-1}g(s,y,z)\mid_{z=-1,s=1}=\lim_{s\to 1^+}(D_s)^N(D_z)^{k-1}g(s,y,z)\mid_{z=-1,s=1}\]
\[=\sum_{n=1}^{N+1}\sum_{j=0}^{k-1}\rho_{n,j,k-1,N+1}(\log\log y)^j(\log y)^{n}+O_{N,k-1}(\exp{(-c''\sqrt{\log y}})).\tag{3.2}\]

\end{theorem}

The key to proving Theorem 3.1 is simply logarithmic differentiation\footnote{Recall that, when a function is differentiable and nonzero, the logarithmic derivative exists, even if the logarithm is undefined, and is independent of the choice of branch cut.} and induction.\footnote{Motivated by the formula for $\Gamma_{N,m}$ in Lemma 2.2, I decided to present a sufficiently elaborate treatment in this section to keep track of information regarding the constants appearing in (3.2) and specifically understand them as members of a field extension of $\mathbb{Q}$ of a finite transcendence degree. } Specifically, for $\sigma>1$, using logarithmic differentiation on (2.1) yields

\[\tag{3.3} \frac{f_z}{f}=\frac{g_z}{g}+\log\zeta.\]

On the other hand, since $f=z^{-1}\prod_{p> y}\bigl(1+\frac{z}{p^s}\bigr)$ we have $\frac{f_z}{f}=-\frac{1}{z}+\sum_{p>y}\frac{1}{p^s+z}$ and substituting above and rearranging gives for $\sigma>1$
\[g_z=(\sum_{p>y}\frac{1}{p^s+z}-\log\zeta(s)-\frac{1}{z})g. \tag{3.4}\]

Applying the Leibniz rule to (3.4) gives for $k-1\ge1$
\begin{align*}\!\!\!
(D_z)^{k-1}g(s,y,z)\,=\,\,&(\sum_{p>y}\frac{1}{p^s+z}-\log\zeta(s)-\frac{1}{z})(D_z)^{k-2}g(s,y,z)\\  
&\!\!\!+\sum_{i=0}^{k-3}{k-2\choose i}\bigg\{ (k-2-i)!(-1)^{k-2-i}(\frac{-1}{z^{k-1-i}}+\sum_{p>y}\frac{1}{(p^s+z)^{k-1-i}})\bigg\} (D_z)^{i}g(s,y,z).\tag{3.5}
\end{align*}

To make way towards (3.2) and (3.3), we must take $s$-derivatives of (3.5). For ease, first rewrite (3.5) using the notation ${(D_z)^{i}g(s,y,z)=\mathcal{G}_i(s,y,z)}$, then substitute $z=-1$ to obtain for $k-1\ge1$
\begin{align*}\hspace*{-0.7 cm}
\!\!\!(D_z)^{k-1}g(s,y,z)\mid_{z=-1}\,=\,\,&\\ 
\!\!\!\mathcal{G}_{k-1}(s,y,-1)=\,\,&(\sum_{p>y}\frac{1}{p^s+(-1)}-\log\zeta(s)-\frac{1}{(-1)})\mathcal{G}_{k-2}(s,y,-1)\\  
&\!\!\!+\sum_{i=0}^{k-3}{k-2\choose i}\bigg\{ (k-2-i)!(-1)^{k-2-i}(\frac{-1}{(-1)^{k-1-i}}+\sum_{p>y}\frac{1}{(p^s+(-1))^{k-1-i}})\bigg\} \mathcal{G}_{i}(s,y,-1)\\
=\,\,&(\sum_{p>y}\frac{1}{p^s-1}-\log\zeta(s)+1)\mathcal{G}_{k-2}(s,y,-1)\\  
&\!\!\!+\sum_{i=0}^{k-3}{k-2\choose i}\bigg\{(k-2-i)!(1+\sum_{p>y} \frac{(-1)^{k-2-i}}{(p^s-1)^{k-1-i}})\bigg\} \mathcal{G}_{i}(s,y,-1).\tag{3.6}
\end{align*}

On the other hand, the series on the right-hand side converges when $s=1$. Thus, we can extend the meaning of the expression on the left-hand side by analytic continuation to a function valid for $\sigma\geq\sigma_0$ with $\sigma_0<1$. In other words, we have an analytic continuation for $\sum_{p>y}\frac{1}{p^s-1}-\log\zeta(s)$ and hence $g_z(s,y,-1)$ for a region $\sigma\geq\sigma_0$ inside the critical strip. By induction on the number of times the $z$-derivative was applied and the relation above, we have for each $i\in\mathbb{N}$ an extension of $\mathcal{G}_i(s,y,-1)$ by analytic continuation to an analytic function on $\sigma\geq\sigma_0$ inside the critical strip. Pick $c$ such that $1-\frac{c}{\log T}\geq\sigma_0$ for $x\geq 1.5$ and so that for $s=\sigma+iT$, $\sigma\geq1-\frac{c}{\log T}$ is in the zero free region and such that all the bounds in \cite[Thm~6.7]{MV} hold. Such a $c$ exists by the classic zero free region and \cite[Thm~6.7]{MV}. Thus, whenever $\sigma\geq\sigma_0$ and hence whenever $\sigma\geq1-\frac{c}{\log T}$ we have for $k-1\ge1$
\begin{align*}
 \!\!\!(D_s)^N\mathcal{G}_{k-1}(s,y,-1)=\,\,&\sum_j{N\choose j}\bigg\{D_s^j(\sum_{p>y}\frac{1}{p^s-1}-\log\zeta(s)+1)(D_s)^{N-j}\mathcal{G}_{k-2}(s,y,-1)\\  
&\!\!\!+\sum_{i=0}^{k-3}{k-2\choose i}(k-2-i)!\big\{D_s^j(1+\sum_{p>y} \frac{(-1)^{k-2-i}}{(p^s-1)^{k-1-i}})\big\}(D_s)^{N-j}\mathcal{G}_{i}(s,y,-1)\bigg\}.\tag{3.7}
\end{align*}

Now first note that $g(s,y,-1)=-\prod_{p\leq y}\frac{1}{1-\frac{1}{p^s}}$ and so 

\[\frac{g_s(s,y,-1)}{g(s,y,-1)}=-\sum_{p\leq y}\frac{\log p}{p^s-1}\Leftrightarrow g_s(s,y,-1)=\bigl(-\sum_{p\leq y}\frac{\log p}{p^s-1}\bigr)g(s,y,-1)\]

and so, by the Leibniz rule  
\[(D_s)^Ng(s,y,-1)=(D_s)^{N-1}g_s(s,y,-1)=\sum_j{N-1\choose j}\bigl(-\sum_{p\leq y}(D_s)^j\frac{\log p}{p^s-1}\bigr)(D_s)^{N-j-1}g(s,y,-1).\]

On the other hand, it can be seen by induction that 
\[-(D_s)^j\frac{\log p}{p^s-1}=(-\log p)^{j+1} \sum_{k=1}^{j+1}\frac{T_{j,k}}{(p^s -1)^k}\]
where $T_{j,k}= k-1!S(j+1,k)$ with $S(j+1,k)$ the Stirling numbers of the second kind (see \cite{TJK}).

This way, we have an explicit recursion for $(D_s)^Ng(s,y,-1)$ in terms of the lower-order derivatives. To get the desired results, we will benefit from upper bounds for $-\sum_{p\leq y}(D_s)^j\frac{\log p}{p^s-1}$ when $\sigma\geq1-\frac{c}{\log T}$ and quantitative estimates for $-\sum_{p\leq y}(D_s)^j\frac{\log p}{p^s-1}\mid_{s=1}$. To do the latter, note that 

\[-\sum_{p\leq y}(D_s)^j\frac{\log p}{p^s-1}\mid_{s=1}=\sum_{p\leq y}\frac{(-\log p)^{j+1}}{p-1}+C_j+\xi_j(y),\]
with $C_j=\sum_p(-\log p)^{j+1} \sum_{k=2}^{j+1}\frac{T_{j,k}}{(p -1)^k}$ and\footnote{An effort to give bounds with explicit constants has been made in several, but not all, cases. Such explicit bounds often merit a dedicated paper; see, e.g., \cite{Ross}.} \[|\xi_j(y)|\leq \sum_{k=2}^{j+1}T_{j,k}\sum_{p>y}\frac{(\log p)^{j+1}}{(p-1)^2}\leq \sum_{k=1}^{j+1}T_{j,k}\sum_{p>y}\frac{(\log p)^{j+1}}{(p-1)^2}\leq25\frac{(j+2!)^2(\frac{\lceil\log y\rceil}{\log 2})^{j+1}}{(j+2)y}\tag{3.8}.\]
Here, the fact that \[1+\sum_{k=2}^{j+1}T_{j,k}=\sum_{k=1}^{j+1}T_{j,k}=2F_{j+1}<2\frac{j+1!}{(\log 2)^{j+2}}\]
follows from the definition of the Fubini numbers and \cite{Fubini}. On the other hand, we also used \[\sum_{p>y}\frac{(\log p)^{j+1}}{(p-1)^2}\leq4\sum_{p>y}\frac{(\log p)^{j+1}}{p^2}\leq 4\int_{y-1}^\infty\frac{(\log u)^{j+1}}{u^2}du\leq 4\frac{j+2!\lceil\log (y)\rceil^{j+1}}{y-1}\leq\frac{4\times19}{9}\frac{j+2!\lceil\log (y)\rceil^{j+1}}{y}\tag{3.9}\]
which follows from an elementary integral estimate.\\

\begin{lemma}
For $N\geq 1$, let $\mathcal{M}_N(y)=\sum_{p\leq y}\frac{(\log p)^N}{p-1}$, then
    \[\mathcal{M}_N(y)=P_{\mathcal{M}_N}(\log y)+O(\exp{(-c'\sqrt{\log y}}))\]
where $P_{\mathcal{M}_N}$ is a polynomial of degree $N$.
\end{lemma}
\begin{proof}
From \cite[Eq. 2.31]{RossSch} we have 
\[\sum_{p\leq y}\frac{\log p}{p}=\log y-B_3 +O(\exp{(-a\sqrt{\log y}}))\]
for some constant $B_3$ whose exact value is in fact well known to equal $B_3=\gamma+\sum_p\frac{\log p}{p(p-1)}$. On the other hand, it is easy to see that 
\[\sum_{p\leq y}\frac{\log p}{p-1}-\sum_{p\leq y}\frac{\log p}{p}=\sum_p\frac{\log p}{p(p-1)} +O(\frac{1}{y}).\]  It directly follows that 
\[\sum_{p\leq y}\frac{\log p}{p-1}=\log y-\gamma +O(\exp{(-a\sqrt{\log y}})).\]
This concludes the base case.

Now, suppose that the claim of the lemma is true up to $N=N'$. Then, by Abel summation or Stieltjes integration,  
\[\mathcal{M}_{N'+1}(y)=\log y \mathcal{M}_{N'}(y)-\int_{1.9}^y\frac{\mathcal{M}_{N'}(t)}{t}dt.\]
The conclusion follows from induction.
\end{proof}
\begin{remark}
Actually, it is already known \cite{Hum} that \[\mathcal{Q}_j(y)=\sum_{p\leq y}\frac{(\log p)^j}{p}=\frac{(\log y)^j}{j}+\beta_j+O(\exp{(-c'\sqrt{\log y}})),\tag{3.10}\]
where the constant $\beta_j$ can be given an explicit formula (see \cite[Prop. 2]{HM} for a related refinement). Note that the same $c'$ and error term can be chosen because, for the constant $\tau_j =\sum_{p}\frac{(\log p)^j}{p(p-1)}$, we have
\[\mathcal{M}_j(y)-\mathcal{Q}_j(y)=\sum_{p\leq y}\frac{(\log p)^j}{p(p-1)}=\tau_j -\sum_{p> y}\frac{(\log p)^j}{p(p-1)},\]
with the smaller error
\[\sum_{p> y}\frac{(\log p)^j}{p(p-1)}\leq\sum_{p>y}\frac{(\log p)^j}{(p-1)^2}\leq \frac{4\times19}{9}\frac{j+1!\lceil\log (y)\rceil^{j}}{y},\]
by (3.9). To connect this with the previous lemma, note that this gives
\[\mathcal{M}_j(y)=\frac{(\log y)^j}{j}+\beta_j+\tau_j+O(\exp{(-c'\sqrt{\log y}})),\tag{3.11}\]
so that we have an explicit formula for $P_{\mathcal{M}_N}(\log y)$ and we see that $P_{\mathcal{M}_N}$ is simply a monomial plus a constant.
\end{remark}

In private communications, K. Alladi has pointed out to me that the error from the quantitative form of the prime number theorem could be used to improve the classic Mertens estimate to 
\[\tag{3.12}\prod_{p\leq y}\frac{1}{1-\frac{1}{p}}=e^\gamma\log y+O(\frac{1}{\exp(c'\sqrt{\log y})}).\]
Actually, \cite[Eq. 2.32]{RossSch} gives 
\[\prod_{p\leq y}\bigl(1-\frac{1}{p}\bigr)=\frac{1}{e^\gamma\log y}+O(\frac{1}{\exp(a\sqrt{\log y})}).\]
On the other hand, if we multiply both sides of the previous equation by $\frac{e^\gamma\log y}{\prod_{p\leq y}(1-\frac{1}{p})} = O((\log y)^2)$, which follows from the classic bounds (see \cite{MV,RossSch}), we deduce (3.12).
\begin{corollary}
 \[(D_s)^Ng(s,y,-1)\mid_{s=1}=\mathcal{P}_{N}(\log y)+O(\exp{(-c'\sqrt{\log y}}))\] 
where $\mathcal{P}_{N}$ is a polynomial of degree $N+1$ with no constant term.
\end{corollary}
\begin{proof}
By induction. The base case is a direct consequence of the refinement of Mertens' estimate stated above. Suppose this is true up to $N'$. We have already noted
\[(D_s)^{N'+1}g(s,y,-1)\mid_{s=1}=(D_s)^{N'}g_s(s,y,-1)\mid_{s=1}=\sum_j{N'\choose j}\bigl(-\sum_{p\leq y}(D_s)^j\frac{\log p}{p^s-1}\mid_{s=1}\bigr)(D_s)^{N'-j}g(s,y,-1)\mid_{s=1}\]
and
\[-\sum_{p\leq y}(D_s)^j\frac{\log p}{p^s-1}\mid_{s=1}=(-1)^{j+1}\mathcal{M}_{j+1}(y)+C_j+\xi_j(y),\]
with $\xi_j(y)=O_j(\frac{(\log y)^{j+1}}{y})$. From the previous lemma,
\[-\sum_{p\leq y}(D_s)^j\frac{\log p}{p^s-1}\mid_{s=1}=(-1)^{j+1}\mathcal{M}'_{j+1}(\log y)+O_j(\exp{(-c'\sqrt{\log y}}))\]
with $\mathcal{M}'_{j+1}$ a polynomial of degree $j+1$. By induction, and the claim just before this corollary, the conclusion follows.
\end{proof}

\begin{remark}
In view of the remark given after Lemma 3.2, and the argument from the previous corollary, the $\mathcal{P}_{N}$ could be given an explicit formula for each $N$ and the coefficients could be seen to belong to a (possibly transcendental) field extension of $\mathbb{Q}$ that can be defined using the constants mentioned so far.
\end{remark}

So far, Corollary 3.3 proves (3.2) in the case $k-1=0$. To complete the proof of (3.2), by induction on $k$, all that remains is to find quantitative estimates for the $s$-derivative of ${ \sum_{p>y}\frac{1}{p^s-1}-\log\zeta(s)+1}$. To begin with, recall that we have seen that there is a $\sigma_0<1$ such that the left-hand side of the following expression 
\[\tag{3.13} \sum_{p>y}\frac{1}{p^s-1}-\log\zeta(s)+1=-\sum_{p\leq y}\frac{1}{p^s}+\sum_{p>y}\frac{1}{p^s(p^s-1)}-\sum_{p}\sum_{m\geq2}\frac{1}{p^{ms}m}+1\]
can be defined by the right-hand side whenever $\sigma\geq\sigma_0$, giving us an analytic function when $\sigma\geq\sigma_0$. Since the above function is analytic for $\sigma\geq\sigma_0$ we have for $j\geq1$
\[\tag{3.14}(D_s)^j(\sum_{p>y}\frac{1}{p^s-1}-\log\zeta(s)+1)\]
\[=-\sum_{p\leq y}\frac{(-\log p)^j}{p^s}+\sum_{p>y}\frac{(-\log p)^j}{p^s(p^s-1)}+\sum_{p>y}(-\log p)^{j} \sum_{k=2}^{j+1}\frac{T_{j,k}}{(p^s -1)^k}-\sum_{p}\sum_{m\geq2}\frac{(-m\log p)^j}{p^{ms}m}.\]

As $\sigma_0<1$ is chosen to guarantee that the above expression is analytic for $\sigma\geq\sigma_0$, it simultaneously guarantees that all the above series converge at $s=1$. In particular, $\sum_{p}\sum_{m\geq2}\frac{(-m\log p)^j}{p^{m}m}=B'_j<\infty$. Now note that 
\[\tag{3.15}\left|\sum_{p>y}\frac{(-\log p)^j}{p^s(p^s-1)}+\sum_{p>y}(-\log p)^{j} \sum_{k=2}^{j+1}\frac{T_{j,k}}{(p^s -1)^k}\right|<\sum_{k=1}^{j+1}T_{j,k}\sum_{p>y}\frac{(\log p)^j}{(p^s-1)^2},\]
and when $s=1$, from (3.8), we have
\[\sum_{k=1}^{j+1}T_{j,k}\sum_{p>y}\frac{(\log p)^j}{(p-1)^2}\leq 25\frac{(j+1!)^2(\frac{\lceil\log y\rceil}{\log 2})^{j}}{(j+1)y}.\]
Putting it all together gives 
 \[(D_s)^j(\sum_{p>y}\frac{1}{p^s-1}-\log\zeta(s)+1)\mid_{s=1}=-\sum_{p\leq y}\frac{(-\log p)^j}{p}-B'_j+O_j(\frac{(\log y)^j}{y}).\]

Combining the above with the quantitative expression for $\sum_{p\leq y}\frac{(\log p)^j}{p}$ stated in the remark after Lemma 3.2 gives for $j\geq 1$
\[(D_s)^j(\sum_{p>y}\frac{1}{p^s-1}-\log\zeta(s)+1)\mid_{s=1}=\mathscr{P}_j(\log y)+O_j(\exp{(-c'\sqrt{\log y}}))\tag{3.16}\]
where $\mathscr{P}_j(\log y)$ is a polynomial in $\log y$ of degree $j$.
\begin{corollary}
Equation (3.2) is true.
\end{corollary}
\begin{proof}
We proceed by strong induction on $k\geq1$. When $k=1 \Leftrightarrow k-1=0$, in view of Corollary 3.3, we can define $\rho_{n,0,0,N+1}$ so that $\mathcal{P}_{N}(\log y)=\sum_{n=1}^{N+1}\rho_{n,0,0,N+1}(\log y)^{n}$. Suppose this is true whenever $0\leq k-1\leq K-2$. 
First, observe that for each $j>0$, $k$ and $k-3\geq i\geq0$ we have  \[D_s^j(1+\sum_{p>y} \frac{(-1)^{k-2-i}}{(p^s-1)^{k-1-i}})\mid_{s=1}\ll_{j,k}\frac{(\log y)^j}{y^{k-1-i}}\leq\frac{(\log y)^j}{y^2} \text{ and }\sum_{p>y} \frac{(-1)^{k-2-i}}{(p^s-1)^{k-1-i}}\mid_{s=1}\ll_{k}\frac{1}{y^{k-1-i}}\leq\frac{1}{y^2}\tag{3.17}.\]
Combining the above with (3.7) gives 
\begin{align*}
 \!\!\!(D_s)^N\mathcal{G}_{K-1}(s,y,-1)\mid_{s=1}=\,\,&\sum_j{N\choose j}D_s^j(\sum_{p>y}\frac{1}{p^s-1}-\log\zeta(s)+1)(D_s)^{N-j}\mathcal{G}_{K-2}(s,y,-1)\mid_{s=1}\\   
&\!\!\!+\sum_{i=0}^{K-3}{K-2\choose i}(K-2-i)! (D_s)^N\mathcal{G}_{i}(1,y,-1)+O_{N,K-1}(\exp{(-c'\sqrt{\log y}})). \tag{3.18}
\end{align*}
By the induction hypothesis, $\sum_{i=0}^{K-3}{K-2\choose i}(K-2-i)! (D_s)^N\mathcal{G}_{i}(1,y,-1)$ already has the desired form. In addition, by induction, we already have the constants $\{\rho_{n,j',K-2,N-j+1}|1\leq n\leq N-j+1,0\leq j'\leq K-2\}$ for $(D_s)^{N-j}\mathcal{G}_{K-2}(s,y,-1)$. Thus, when $j\geq 1$, combining the preceding inductive hypothesis with (3.16) yields 
\begin{align*}
 \!\!\!&(D_s)^j(\sum_{p>y}\frac{1}{p^s-1}-\log\zeta(s)+1)(D_s)^{N-j}\mathcal{G}_{K-2}(s,y,-1)\mid_{s=1}\,\,\\  
 &\quad=\mathscr{P}_j(\log y)\sum_{n=1}^{N-j+1}\sum_{j'=0}^{K-2}\rho_{n,j',K-2,N-j+1}(\log\log y)^{j'}(\log y)^{n}+O_{N,K-1}(\exp{(-c''\sqrt{\log y}})),\tag{3.19}
\end{align*}
where $c''$ is chosen to be smaller than $c'$.

For $j=0$, note that, from (3.13) and the well-known bounds on $\sum_{p\leq y}\frac{1}{p}$ (see, e.g., \cite[Eq. 2.30]{RossSch}), we have
 \[\lim_{s\to1^+}\sum_{p>y}\frac{1}{p^s-1}-\log\zeta(s)+1=1-\gamma-\log\log y +O(\exp{(-a\sqrt{\log y}})).\tag{3.20}\]
Thus, the inductive hypothesis and (3.20) directly give
\begin{align*}
 \!\!\!&(\sum_{p>y}\frac{1}{p^s-1}-\log\zeta(s)+1)(D_s)^{N}\mathcal{G}_{K-2}(s,y,-1)\mid_{s=1}\,\,\\  
 &\quad=(1-\gamma-\log\log y)\sum_{n=1}^{N+1}\sum_{j'=0}^{K-2}\rho_{n,j',K-2,N+1}(\log\log y)^{j'}(\log y)^{n}+O_{N,K-1}(\exp{(-c''\sqrt{\log y}})).\tag{3.21}
\end{align*}

The conclusion follows by substituting the results of (3.19) and (3.21) back into (3.18).
\end{proof}

\begin{remark}
As was previously mentioned, given the recursive nature of (3.18), the $\rho_{n,j,k-1,N+1}$ can be computed explicitly and likewise could be seen to belong to a (possibly transcendental) field extension of $\mathbb{Q}$ that can be defined using the constants mentioned so far. The explicit formulas are quite complicated, but I will point out a striking feature of $\rho_{1,j,k-1,1}$, the case where no $s$-derivatives are taken, which relates to the case $N=1$ in (1.2). Note that if there are no 
$s$-derivatives, then the first recursive relation in the proof becomes
\[ \mathcal{G}_{k-1}(1,y,-1)=(1-\gamma-\log\log y)\mathcal{G}_{k-2}(1,y,-1)+\sum_{i=0}^{k-3}{k-2\choose i}(k-2-i)! \mathcal{G}_{i}(1,y,-1)+O_{N,k-1}(\exp{(-c''\sqrt{\log y}})). \]
Because $\mathcal{G}_{0}(1,y,-1)=-e^\gamma \log y + O(\frac{1}{\exp(c'\sqrt{\log y})})$, and the factor $(1-\gamma-\log\log y)$ extends the transcendence of the field by at most $\gamma$, a straightforward induction shows that $\rho_{1,j,k-1,1}\in\mathbb{Q}(\gamma,e^{\gamma})$ for all $k-1,j\geq0$. Remember, the $k$ here represents a duality with the $k$-th prime factor, so what is seen here is that for $k\geq 3$, the contribution to the transcendence from $\mathcal{G}_{k-1}(1,y,-1)$ is redundant. We will return to this after proving the main theorem.
\end{remark}

We now turn our attention to equation (3.1).

\begin{lemma}
    There is a constant $\mathscr{C}$ (that may depend on $Y_0$, $\mathscr{p}$, and $T_0$ but is otherwise absolute) such that for $\sigma\geq1-\frac{c}{\log T}$, $y\leq Y_0T^\mathscr{p}$, and $T\geq T_0$ we have 
    \[\sum_{p\leq y}\frac{1}{p^\sigma}\leq\log\log y +\mathscr{C} .\]

\end{lemma}

\begin{proof}
I will follow the approach in \cite{MV}, but give more explicit estimates for the specific domain in question. First, without loss of generality, $\sigma\leq 1$. Now, just like the proof of \cite[Lem. 7.3]{MV}

\[\sum_{p\leq y}\frac{1}{p^\sigma}=\int_{1.9}^{y}\frac{1}{u^{\sigma}}d\pi(u)=\int_{1.9}^y\frac{1}{u^\sigma\log u}du+\int_{1.9}^y\frac{1}{u^\sigma}dr(u).\]

For the first integral, note that for $1\leq u\leq Y_0 T^\mathscr{p}$ we have the inequality \[{u^{\frac{c}{\log T}}\leq1+(Y_0^{\frac{c}{\log T}}e^{c\mathscr{p}}-1)\frac{\log u}{\log Y_0 T^\mathscr{p}}}<1+Y_0^{\frac{c}{\log T_0}}e^{c\mathscr{p}}\frac{\log u}{\log Y_0 T^\mathscr{p}}.\]

\noindent Thus, when $y\leq Y_0 T^\mathscr{p}$
\[\int_{1.9}^y\frac{du}{u^\sigma\log u}<\int_{1.9}^y\frac{du}{u\log u}+\int_{1.9}^{Y_0 T^\mathscr{p}}\frac{Y_0^{\frac{c}{\log T_0}}e^{c\mathscr{p}}\frac{\log u}{\log Y_0 T^\mathscr{p}}}{u\log u}du<\log\log y+Y_0^{\frac{c}{\log T_0}}e^{c\mathscr{p}}.\]

For the second integral, we start by integrating by parts to get
\[\int_{1.9}^y\frac{1}{u^\sigma}dr(u)=\frac{r(y)}{y^\sigma}-\frac{r(1.9)}{1.9^\sigma}+\sigma\int_{1.9}^y\frac{r(u)}{u^{\sigma+1}}du.\]
Now, note that
\[-\frac{r(1.9)}{1.9^\sigma}=\frac{\mathrm{li}(1.9)}{1.9^\sigma}<\mathrm{li}(1.9);\]
in addition, recall that the quantitative form of the prime number theorem gives constants $A$ and $A'$ such that $|r(y)|\leq\frac{Ay}{\exp({A'\sqrt{\log y}})}$. It follows that, for $y\geq1.9$, we have
\[\left|\frac{r(y)}{y^\sigma}\right|\leq A\exp{(\frac{c\log y}{\log T}-A'\sqrt{\log y})}\leq A \frac{Y_0^{\frac{c}{\log T_0}}e^{c\mathscr{p}}}{\exp{(A'\sqrt{\log y})}}\leq A\frac{Y_0^{\frac{c}{\log T_0}}e^{c\mathscr{p}}}{\exp{(A'\sqrt{\log 1.9})}}.\]Lastly, as $|r(y)|\leq\frac{Ay}{\exp{A'\sqrt{\log y}}}$, by the trinagle inequality for integrals
\[\sigma\left|\int_{1.9}^y\frac{r(u)}{u^{\sigma+1}}du\right|\leq\sigma\int_{1.9}^y\frac{A}{u^\sigma\exp{(A'\sqrt{\log u})}}du\leq A\int_{1.9}^y\frac{Y_0^{\frac{c}{\log T_0}}e^{c\mathscr{p}}}{u\exp{(A'\sqrt{\log u}})}du\leq \mathscr{C}'',\]where $\mathscr{C}''$ is an absolute constant (sub $v=\sqrt{\log u}$). Then $\mathscr{C}$ is the sum of the three constants.

\end{proof}

\begin{corollary}
Let $\mathscr{C}$ be the constant from the previous lemma. For $\sigma\geq1-\frac{c}{\log T}$, $y\leq Y_0T^\mathscr{p}$, and $T\geq T_0$ we have  \[|g(s,y,-1)|\leq e^{4\zeta(2.1\sigma_0)+\mathscr{C}}\log y. \]
\end{corollary}
\begin{proof}
 \[|g(s,y,-1)|\leq e^{4\zeta(2.1\sigma_0)}\exp(\sum_{p\leq y}\frac{1}{p^\sigma})\ \leq e^{4\zeta(2.1\sigma_0)+\mathscr{C}}\log y\]
\end{proof}
\begin{lemma}
There are constants $\mathscr{c}_j$ such that for $\sigma\geq1-\frac{c}{\log T}$, $1.9\leq y\leq Y_0T^\mathscr{p}$, and $T\geq T_0$: if $j\geq 1$
   \[ \left|(D_s)^j(\sum_{p>y}\frac{1}{p^s-1}-\log\zeta(s)+1)\right|\leq\frac{3Y_0^{\frac{c}{\log T_0}}e^{c\mathscr{p}}}{\log 1.5}(\log y)^{j}+\mathscr{c}_j,\]
otherwise
\[ \left|\sum_{p>y}\frac{1}{p^s-1}-\log\zeta(s)+1\right|\leq\log \log y+\mathscr{c}_0.\]
\end{lemma}
\begin{proof}
In the case $j=0$ by (3.13) and Lemma 3.5 we have 
\[\left|\sum_{p>y}\frac{1}{p^s-1}-\log\zeta(s)+1\right|\leq \sum_{p\leq y}\frac{1}{p^\sigma}+\sum_{p>y}\frac{1}{p^{\sigma_0}(p^{\sigma_0}-1)}+\sum_{p}\sum_{m\geq2}\frac{1}{p^{m\sigma_0}m}\]
\[\leq \log\log y +\mathscr{C} +\sum_{p>y}\frac{1}{p^{\sigma_0}(p^{\sigma_0}-1)}+\sum_{p}\sum_{m\geq2}\frac{1}{p^{m\sigma_0}m}. \]
From (3.14) , (3.15), and the triangle inequality we have for $j>0$
      \[\left|(D_s)^j(\sum_{p>y}\frac{1}{p^s-1}-\log\zeta(s)+1)\right|\leq\sum_{p\leq y}\frac{(\log p)^j}{p^\sigma}+\sum_{k=1}^{j+1}T_{j,k}\sum_p\frac{(\log p)^j}{(p^{\sigma_0}-1)^2}+\sum_{p}\sum_{m\geq2}\frac{(m\log p)^j}{p^{m\sigma_0}m}\]

      \[\leq Y_0^{\frac{c}{\log T_0}}e^{c\mathscr{p}}(\log y)^{j-1}\sum_{p\leq y}\frac{\log p}{p}+\sum_{k=1}^{j+1}T_{j,k}\sum_p\frac{(\log p)^j}{(p^{\sigma_0}-1)^2}+\sum_{p}\sum_{m\geq2}\frac{(m\log p)^j}{p^{m\sigma_0}m}.\]
      On the other hand, it is well known that $|\sum_{p\leq y}\frac{\log p}{p}-\log y|\leq 2$ for $y\geq 1$. Thus, for $y\geq1.9$
      \[(\log y)^{j-1}\sum_{p\leq y}\frac{\log p}{p}\leq (\log y)^{j-1}(2+\log y )\leq (\log y)^{j-1}(2\frac{\log y}{\log 1.5}+\log y )\leq \frac{3}{\log 1.5}(\log y)^{j}.\]
      The conclusion follows.
\end{proof}
\begin{lemma}
For $j\geq0$, there are constants $\mathscr{c}'_j$ such that for $\sigma\geq1-\frac{c}{\log T}$, $1.9\leq y\leq Y_0T^\mathscr{p}$, and $T\geq T_0$
\[\left|-\sum_{p\leq y}(D_s)^j\frac{\log p}{p^s-1}\right |\leq\frac{3Y_0^{\frac{c}{\log T_0}}e^{c\mathscr{p}}}{\log 1.5}(\log y)^{j+1}+\mathscr{c}'_j.\]
\end{lemma}
\begin{proof}
The fact that \(-(D_s)^j\frac{\log p}{p^s-1}=(-\log p)^{j+1} \sum_{k=1}^{j+1}\frac{T_{j,k}}{(p^s -1)^k}\) was already noted in the discussion preceding Lemma 3.2. Hence,
\[\left|-\sum_{p\leq y}(D_s)^j\frac{\log p}{p^s-1}\right |\leq\sum_{p\leq y}\frac{(\log p)^{j+1}}{p^\sigma-1}+\sum_{p}(\log p)^{j+1} \sum_{k=2}^{j+1}\frac{T_{j,k}}{(p^{\sigma_0} -1)^2}\]
\[\leq\sum_{p\leq y}\frac{(\log p)^{j+1}}{p^\sigma}+\sum_{p}\frac{(\log p)^{j+1}}{(p^{\sigma_0}-1)p^{\sigma_0}}+\sum_{p}(\log p)^{j+1} \sum_{k=2}^{j+1}\frac{T_{j,k}}{(p^{\sigma_0} -1)^2}\]
\[\leq\frac{3Y_0^{\frac{c}{\log T_0}}e^{c\mathscr{p}}}{\log 1.5}(\log y)^{j+1}+\sum_{p}\frac{(\log p)^{j+1}}{(p^{\sigma_0}-1)p^{\sigma_0}}+\sum_{p}(\log p)^{j+1} \sum_{k=2}^{j+1}\frac{T_{j,k}}{(p^{\sigma_0} -1)^2}.\]
\end{proof}
\begin{corollary}
Inequality (3.1) is true.    
\end{corollary}
\begin{proof}
    We now proceed by strong induction on $N+k-1$. I will look at two cases, one where no $z$-derivatives are taken, and one where there is at least one occurrence of a $z$-derivative. If there are no $z$-derivatives, we proceed by induction. Suppose this is true up to $N'$. Because \[(D_s)^{N'+1}g(s,y,-1)=\sum_j{N'\choose j}\bigl(-\sum_{p\leq y}(D_s)^j\frac{\log p}{p^s-1}\bigr)(D_s)^{N'-j}g(s,y,-1),\]
    the previous lemma along with the induction hypothesis give
    \[\left|(D_s)^{N'+1}g(s,y,-1)\right|\leq\sum_j{N'\choose j}\left|-\sum_{p\leq y}(D_s)^j\frac{\log p}{p^s-1}\right| \left|(D_s)^{N'-j}g(s,y,-1)\right|\]
    \[\leq\sum_j{N'\choose j}\left(\frac{3Y_0^{\frac{c}{\log T_0}}e^{c\mathscr{p}}}{\log 1.5}(\log y)^{j+1}+\mathscr{c}'_j\right) m_{N',0}(\log y)^{N'-j+1}\]
    \[\leq m_{N',0}(\log y)^{N'+2}\sum_j{N'\choose j}\left(\frac{3Y_0^{\frac{c}{\log T_0}}e^{c\mathscr{p}}}{\log 1.5}+\frac{\mathscr{c}'_j}{(\log 1.5)^{j+1}}\right).\]
    Hence we can choose \[m_{N'+1,0}= m_{N',0}(1+\sum_j{N'\choose j}\left(\frac{3Y_0^{\frac{c}{\log T_0}}e^{c\mathscr{p}}}{\log 1.5}+\frac{\mathscr{c}'_j}{(\log 1.5)^{j+1}}\right)).\]
    This settles the case with no $z$-derivatives.\\ 
    
    We now address the case where there is at least one occurrence of a $z$-derivative. Suppose (3.1) is true up to $N'+k'-2$. Because there is at least one application of a $z$-derivative, we can invoke (3.7) to get

\begin{align*}\hspace*{-1.1 cm}
\!\!\!\!\!\left|(D_s)^{N'}(D_z)^{k'-1}g(s,y,z)\mid_{z=-1}\right|\,=\,\,&\\ 
\!\!\!\!\!\left|(D_s)^{N'}\mathcal{G}_{k'-1}(s,y,-1)\right|=\,\,&\!\bigg|\sum_j{N'\choose j}\bigg\{D_s^j(\sum_{p>y}\frac{1}{p^s-1}-\log\zeta(s)+1)(D_s)^{N'-j}\mathcal{G}_{k'-2}(s,y,-1)\\  
&\!\!\!+\sum_{i=0}^{k'-3}{k'-2\choose i}(k'-2-i)!\big\{D_s^j(1+\sum_{p>y} \frac{(-1)^{k'-2-i}}{(p^s-1)^{k'-1-i}})\big\}(D_s)^{N'-j}\mathcal{G}_{i}(s,y,-1)\bigg\}\bigg|\\
\!\!\!\!\leq\,&\sum_j{N'\choose j}\bigg\{\mathscr{r}_j'+\bigg|D_s^j(\sum_{p>y}\frac{1}{p^s-1}-\log\zeta(s)+1)\bigg|\bigg\}\\
&\times m_{N',k'-2}(\log y)^{N'-j+1}(\log\log (y+1))^{k'-2}\\  
\leq&\, m'_{N',k'-1}(\log y)^{N'+1}(\log\log (y+1))^{k'-1}.
\end{align*}
The last line follows from Lemma 3.7. In the above, for each $j\leq N'$, $\mathscr{r'}_{j}$ are some constants that depend only on $N'$ and $k'-1$ whose existence is guaranteed by the induction hypothesis.\footnote{To write down a formula for $\mathscr{r'}_{j}$ one can start by bounding $|D_s^j(1+\sum_{p>y} \frac{(-1)^{k'-2-i}}{(p^s-1)^{k'-1-i}})|$ in the relevant region and combining this with the bounds of $(D_s)^{j}\mathcal{G}_{j'}(s,y,-1)$ assumed by induction.} Note that in the instance of $j=0$ we can bound $\log\log y +\mathscr{c}_0$ from above by $(100|\mathscr{c}_0|+1)\log\log (y+1)$. The conclusion follows by letting ${m_{N',k'-1}=\max{\{m'_{N',k'-1},m_{N'-1,k'-1},m_{N',k'-2}\}}}$.

\end{proof}

\section{The quantitative statement for $y\leq Y_0\exp(\mathscr{p}\frac{\log x}{(\log\log (x+1))^{1+\epsilon}})$}

At this point, I will adopt the more standard approach (see \cite[p. 230]{MV} and also \cite{KA82}). Let $b=1-\frac{c}{\log T}$ and\footnote{Note that $T\geq T_0=\exp{(\frac{c}{\ell(N)}\frac{\log 1.9}{\log \log 2.9})}$. } $T=\exp{(\frac{c}{\ell(N)}\frac{\log x}{\log \log (x+1)})}$ with $\ell(N)=N+7$. Replace the contour in (2.3) with
\begin{enumerate}[(i)]
    \item horizontal lines connecting $a\pm iT$ to $b\pm iT$,
    \item vertical lines connecting $b\pm iT$ to $b\pm \frac{i}{\log{x}}$, and
    \item the important contour composed of horizontal lines connecting $b\pm \frac{i}{\log{x}}$ to $1\pm \frac{i}{\log{x}}$ and the semi circle $\{1+ \frac{e^{i\theta}}{\log{x}}:-\frac{\pi}{2}\leq\theta\leq\frac{\pi}{2}\}$.
\end{enumerate}

For the remainder of this paper, we will come back to the following simple relation obtained from Leibniz's rule. 
\[\frac{(-1)^k}{k-1!}(D_z)^{k-1}f(s,y,z)\mid_{z=-1}=\frac{(-1)^k}{k-1!\zeta(s)}\sum_{i=0}^{k-1}{k-1\choose i}(\log\zeta(s))^{k-i-1}(D_z)^{i}g(s,y,z)\mid_{z=-1}\tag{4.1}\]

\begin{lemma}
Suppose $y\leq Y_0\exp{(\mathscr{p}\frac{\log x}{\log\log (x+1)})}$, the combined contribution from (i) and (ii) is controlled by the error term in Theorem 1.1.
\end{lemma}
\begin{proof}
By replacing $\mathscr{p}$ with $\frac{\mathscr{p}\ell(N)}{c}$ if necessary, we invoke Theorem 3.1 to get
\[\tag{4.2}\left|\frac{(-1)^k}{k-1!}(D_z)^{k-1}f(s,y,z)\mid_{z=-1}\right|\leq\frac{\log y(\log\log (y+1))^{k-1}}{|\zeta(s)|}\sum_{i=0}^{k-1}(\log\zeta(s))^{k-i-1}m_{0,k-1}\]
\[\leq km_{0,k-1}\frac{\log y(\log\log (y+1))^{k-1}}{|\zeta(s)|}\lceil|\log\zeta(s)|\rceil^{k-1}.\]

On the other hand, by choice of $c$ and \cite[Thm 6.7]{MV} there is a constant $C_{mv}$ such $|\frac{1}{\zeta(s)}|\leq C_{mv}\log T$ and $\lceil|\log\zeta(s)|\rceil\leq C_{mv} \log\log T$ everywhere in the shifted contour. With this in mind, and the fact that $y\leq x$, we have 
\[\tag{4.3}\left|\frac{(-1)^k}{k-1!}(D_z)^{k-1}f(s,y,z)\mid_{z=-1}\right|\leq kC_{mv}^km_{0,k-1}\log x(\log\log (x+1))^{k-1} \log T(\log\log T)^{k-1}\]
everywhere in the shifted contour.

Thus, the contribution from (i) is controlled because
\[ \left|\frac{1}{2\pi i}\int_{a\pm iT}^{b\pm iT}\frac{(-1)^k}{k-1!}(D_z)^{k-1}f(s,y,z)\mid_{z=-1}\frac{x^s}{s}ds\right|\]
\[\leq kC_{mv}^km_{0,k-1}|b-a|\frac{x^a}{T} \log x(\log\log (x+1))^{k-1} \log T(\log\log T)^{k-1}\]
\[\ll_k\frac{x}{\exp(\frac{c}{2\ell(N)}\frac{\log x}{\log\log (x+1)})}.\]
In a similar manner, the contribution from (ii) is controlled because
\[\hspace*{-1cm} \left|\frac{1}{2\pi i}\int_{b\pm \frac{i}{\log x}}^{b\pm iT}\frac{(-1)^k}{k-1!}(D_z)^{k-1}f(s,y,z)\mid_{z=-1}\frac{x^s}{s}ds\right|\]
\[ \leq kC_{mv}^km_{0,k-1} \log x(\log\log (x+1))^{k-1} \log T(\log\log T)^{k-1}x^b\int_\frac{1}{\log x}^T \frac{1}{t} dt\]
\[\leq kC_{mv}^km_{0,k-1} \log x(\log\log (x+1))^{k-1} \log T(\log\log T)^{k-1}\frac{x}{(\log x)^{\ell(N)}}\int_\frac{1}{\log x}^T \frac{1}{t} dt\]
\[\ll_k\frac{x(\log\log x)^{2k-2}}{(\log x)^{\ell(N)-3}}=\frac{x(\log\log x)^{2k-2}}{(\log x)^{N+4}}\ll_k\frac{x}{(\log x)^{N+3.5}}.\]
\end{proof}

 At last, we can now turn our attention to the main contribution, which comes from the important contour $\mathcal{C}(x)$. Recall that I denote
$(D_z)^{i}g(s,y,z)=\mathcal{G}_i(s,y,z)$. We still need to be careful and have some foresight. In particular, we must expand a bit beyond the desired order of the quantitative statement. Specifically, let\footnote{This choice of $N'=N+\lceil\frac{N+3+k}{\epsilon}\rceil$ is important to keep in mind for Lemma 4.3.} $N'=N+\lceil\frac{N+3+k}{\epsilon}\rceil$ and use Taylor's formula (see, e.g., \cite{Rudin}) to get
\[\tag{4.4}\mathcal{G}_i(s,y,-1)=\sum_{j=0}^{N'}\frac{(D_s^j\mathcal{G}_i)(1,y,-1)}{j!}(s-1)^j+\frac{1}{N'!}\int_1^s(t-1)^{N'}(D_t^{N'+1})\mathcal{G}_i(t,y,-1)dt,\]
where the notational shorthand $(D_s^j\mathcal{G}_i)(1,y,-1)=(D_s^j)\mathcal{G}_i(s,y,-1)|_{s=1}$ is used. Substituting the previous equation in equation (4.1) gives 
\begin{align*}\!\!\!
\frac{(-1)^k}{k-1!s}(D_z)^{k-1}f(s,y,z)\mid_{z=-1}=&\frac{(-1)^k}{k-1!s\zeta(s)}\sum_{i=0}^{k-1}{k-1\choose i}(\log\zeta(s))^{k-i-1}(\sum_{j=0}^{N'}\frac{(D_s^j\mathcal{G}_i)(1,y,-1)}{j!}(s-1)^j)\\  
&\!\!\!+\frac{(-1)^k}{k-1!s\zeta(s)}\sum_{i=0}^{k-1}{k-1\choose i}\frac{(\log\zeta(s))^{k-i-1}}{N'!}\int_1^s(t-1)^{N'}(D_t^{N'+1})\mathcal{G}_i(t,y,-1)dt \tag{4.5}
.\end{align*}

Let 
\[\tag{4.6}\varepsilon_1(s,y)=\frac{(-1)^k}{k-1!s\zeta(s)}\sum_{i=0}^{k-1}{k-1\choose i}\frac{(\log\zeta(s))^{k-i-1}}{N'!}\int_1^s(t-1)^{N'}(D_t^{N'+1})\mathcal{G}_i(t,y,-1)dt.\]

The crucial role of controlling the contribution from (4.6) in the important contour $\mathcal{C}(x)$ is played by the following lemma.

\begin{lemma}
There is a constant $F'_{N,k-1}$ that does not depend on any parameters other than $N$, $k-1$, $Y_0$, $\epsilon$ and $\mathscr{p}$, such that \[\left|\frac{1}{2\pi i}\int_{\mathcal{C}(x)}\varepsilon_1(s,y)x^sds\right|\leq F'_{N,k-1}\frac{x}{\log x}(\frac{\log y}{\log x})^{N+2} (\log\log (x+1))^{k-1}.\]    
\end{lemma}
\begin{proof}
Expanding the (4.6) using $\log\zeta(s)=\log\frac{1}{s-1}+\Delta(s)$, which is valid in the convex hull of $\mathcal{C}(x)$, gives 
\[\frac{(-1)^k}{k-1!s\zeta(s)}\sum_{i=0}^{k-1}{k-1\choose i}\frac{(\log\zeta(s))^{k-i-1}}{N'!}\int_1^s(t-1)^{N'}(D_t^{N'+1})\mathcal{G}_i(t,y,-1)dt=\]
\[\frac{(-1)^k}{k-1!N'!s\zeta(s)}\sum_{i,r\geq 0}{k-1\choose i,r,k-i-r-1}(\Delta(s))^{k-i-r-1}(\log\frac{1}{s-1})^r\int_1^s(t-1)^{N'}(D_t^{N'+1})\mathcal{G}_i(t,y,-1)dt.\]
I will show that the contribution from the second sum is controlled by the desired error by looking at its terms one by one. Considering
\[\left|\frac{1}{N'!i!r!k-i-r-1!}\int_{\mathcal{C}(x)}(\Delta(s))^{k-i-r-1}(\log\frac{1}{s-1})^r\int_1^s(t-1)^{N'}(D_t^{N'+1})\mathcal{G}_i(t,y,-1)dt\frac{x^s}{s\zeta(s)}ds\right|\]
and noting that $|\frac{1}{s\zeta(s)}|\leq M_\zeta|1-s|$ everywhere in the important contour, the above is
\[\hspace*{-1cm}\leq e x(\frac{2c}{\log T}+\frac{\pi}{2\log x}) M_{\zeta}(M_\Delta)^{k-i-r-1}(\log\log(x+1)+\ell(N)+2\pi)^r\max_{s\in\mathcal{C}(x)}\{|s-1|^{N'+2}\}\max_{s\in\mathcal{C}(x)}\{|D_t^{N'+1}\mathcal{G}_i(t,y,-1)\mid_{t=s}|\}\]
\[\hspace*{-1cm}\leq e x(\frac{2c}{\log T}+\frac{\pi}{2\log x}) M_{\zeta}(M_\Delta)^{k-i-r-1}(\log\log(x+1)+\ell(N)+2\pi)^r(\frac{c}{\log T}+\frac{1}{\log x})^{N'+2}\max_{s\in\mathcal{C}(x)}\{|D_t^{N'+1}\mathcal{G}_i(t,y,-1)\mid_{t=s}|\}.\]

Because $y\leq Y_0\exp(\mathscr{p}\frac{\log x}{|\log\log (x+1)|^{1+\epsilon}})$, the previous section guarantees a constant $m_{N'+1,i}$ the depends only on $i$ and $N'+1$ such that the above is 
\[\leq ex(\frac{2c}{\log T}+\frac{\pi}{2\log x})m_{N'+1,i}M_\zeta (M_\Delta)^{k-i-r-1}(\log\log(x+1)+\ell(N)+2\pi)^r(\frac{1+\ell(N)\log\log (x+1)}{\log x})^{N'+2}\]
\[\times(\log y)^{N'+2}(\log\log (y+1))^{i}\]
\[\hspace*{-1cm}\leq6 x m_{N'+1,i}M_\zeta (M_\Delta)^{k-i-r-1}(\log\log(x+1)+\ell(N)+2\pi)^r(\frac{1+\ell(N)\log\log (x+1)}{\log x})^{N'+3}(\log y)^{N'+2}(\log\log (y+1))^{i}. \tag{4.7}\]

Because $y\leq Y_0\exp(\mathscr{p}\frac{\log x}{|\log\log (x+1)|^{1+\epsilon}})$ and $x,y\geq1.9$ we can find $Y_0'$ such that 
\[\log y \le \log Y_0 +\mathscr{p}\frac{\log x}{|\log\log (x+1)|^{1+\epsilon}}\le Y_0'\mathscr{p}\frac{\log x}{|\log\log (x+1)|^{1+\epsilon}}\]
whenever $1.9\le y\le Y_0\exp(\mathscr{p}\frac{\log x}{|\log\log (x+1)|^{1+\epsilon}})$.
The choice $N'=N+\lceil\frac{N+3+k}{\epsilon}\rceil$ now plays a crucial role because (4.7) is
\begin{align*}
 &\leq(Y_0'\mathscr{p})^{\lceil\frac{N+3+k}{\epsilon}\rceil}6x m_{N'+1,i}M_\zeta (M_\Delta)^{k-i-r-1}(\log\log(x+1)+\ell(N)+2\pi)^r\\
 &\quad\times(\frac{1+\ell(N)\log\log (x+1)}{\log x})^{\lceil\frac{N+3+k}{\epsilon}\rceil+N+3}(\log y)^{N+2}(\log\log (y+1))^{i} (\frac{(\log x)^{\lceil\frac{N+3+k}{\epsilon}\rceil}}{(\log \log (x+1))^{\lceil\frac{N+3+k}{\epsilon}\rceil (1+\epsilon)}})\\
  &\leq(Y_0'\mathscr{p})^{\lceil\frac{N+3+k}{\epsilon}\rceil}6 m_{N'+1,i}M_\zeta (M_\Delta)^{k-i-r-1} \frac{x}{\log x}(\frac{\log y}{\log x})^{N+2}(\log\log (y+1))^{i} \\
 &\quad\times(\log\log(x+1)+\ell(N)+2\pi)^r\frac{({1+\ell(N)\log\log (x+1)})^{\lceil\frac{N+3+k}{\epsilon}\rceil+N+3}}{(\log \log (x+1))^{\lceil\frac{N+3+k}{\epsilon}\rceil (1+\epsilon)}}\\
  &\leq(Y_0'\mathscr{p})^{\lceil\frac{N+3+k}{\epsilon}\rceil}6 m_{N'+1,i}M_\zeta (M_\Delta)^{k-i-r-1} \frac{x}{\log x}(\frac{\log y}{\log x})^{N+2}(\log\log (x+1))^{i+r}\\
 &\quad\times \frac{(\log\log(x+1)+\ell(N)+2\pi)^r}{(\log\log (x+1))^r}\frac{({1+\ell(N)\log\log (x+1)})^{\lceil\frac{N+3+k}{\epsilon}\rceil+N+3}}{(\log \log (x+1))^{\lceil\frac{N+3+k}{\epsilon}\rceil (1+\epsilon)}}.
\end{align*}

The lemma follows from the fact that $i+r\leq k-1$ and \[ \frac{(\log\log(x+1)+\ell(N)+2\pi)^r}{(\log\log (x+1))^r}\frac{({1+\ell(N)\log\log (x+1)})^{\lceil\frac{N+3+k}{\epsilon}\rceil+N+3}}{(\log \log (x+1))^{\lceil\frac{N+3+k}{\epsilon}\rceil (1+\epsilon)}}\ll1.\]
\end{proof}

The hard part is already behind us, as the remaining contributions are much easier to control. Since we have already accounted for (4.6) we need only to account for the first sum in (4.5). We can write $\log \zeta(s)=\log\frac{1}{s-1}+\Delta(s)$ with $\Delta(s)$ analytic in the relevant region and substitute to get for any $n\geq0$ 
\[\frac{(-1)^k}{k-1!s\zeta(s)}\sum_{0\leq i,r\leq k-1}{k-1\choose i,r,k-i-r-1}(\Delta(s))^{k-i-r-1}(\log\frac{1}{s-1})^r(\sum_{j=0}^{N'}\frac{(D_s^j\mathcal{G}_i)(1,y,-1)}{j!}(s-1)^j)\]
\[=\sum_{0\leq i,r\leq k-1}(\Theta_{i,r,k-1,n}(s)+E_{i,r,k-1,n}(s))(\log\frac{1}{s-1})^r(\sum_{j=0}^{N'}\frac{(D_s^j\mathcal{G}_i)(1,y,-1)}{j!}(s-1)^j),\]
where
\[\Theta_{i,r,k-1,n}(s)=\frac{(-1)^k{k-1\choose i,r,k-i-r-1}}{k-1!}\left(\sum_{j=0}^{n}\frac{(D_s)^j\frac{1}{s\zeta(s)}|_{s=1}}{j!}(s-1)^j\right)\left(\sum_{j=0}^{n}\frac{(D_s^j\Delta)(1)}{j!}(s-1)^j\right)^{k-i-r-1}\tag{4.8},\]
and $E_{i,r,k-1,n}(s)=\frac{(-1)^k{k-1\choose i,r,k-i-r-1}}{k-1!}\frac{(\Delta(s))^{k-i-r-1}}{s\zeta(s)}-\Theta_{i,r,k-1,n}(s)$. Let $\phi_{i,J,n}(s)$ be such that \par\pagebreak[2]
\[\sum_{0\leq i,r\leq k-1}\Theta_{i,r,k-1,n}(s)(\log\frac{1}{s-1})^r(\sum_{j=0}^{N'}\frac{(D_s^j\mathcal{G}_i)(1,y,-1)}{j!}(s-1)^j)=\sum_{\substack{1\leq i\leq I_{n,N',k-1}\\0\leq J\leq k-1} }\phi_{i,J,n}(y)(s-1)^i(\log{\frac{1}{s-1}})^J\tag{4.9},\] 
with $I_{n,N',k-1}$ depending only on $n$, $N'$ and $k-1$. Such an $I_{n,N',k-1}$ exists because each $\Theta_{i,r,k-1,n}(s)$ is a polynomial and the sum in the LHS of (4.9) is finite.
It follows that
\[\frac{(-1)^k}{k-1!s}(D_z)^{k-1}f(s,y,z)\mid_{z=-1}=\sum_{\substack{ 1\le i\le I\\0\leq J\leq k-1} }\phi_{i,J,N+10}(y)(s-1)^i(\log{\frac{1}{s-1}})^J+\varepsilon_2(s,y),\]

where $I=I_{N+10,N',k-1}$ and $\varepsilon_2(s,y)=\varepsilon_1(s,y)+\varepsilon'(s,y)$ with
\[\varepsilon'(s,y)=\sum_{0\leq i,r\leq k-1}E_{i,r,k-1,N+10}(s)(\log\frac{1}{s-1})^r(\sum_{j=0}^{N'}\frac{(D_s^j\mathcal{G}_i)(1,y,-1)}{j!}(s-1)^j).\]

Again, because $\varepsilon_1(s,y)$ has already been accounted for, we need only to account for the other terms. This is done by the following lemma. 
\begin{lemma}
    There is a pair of constants, $F''_{N,k-1}$ and $F'''_{N,k-1}$, that do not depend on any parameters other than $N$, $k-1$, $Y_0$, $\epsilon$ and $\mathscr{p}$, such that \[\left|\frac{1}{2\pi i}\int_{\mathcal{C}(x)}\varepsilon'(s,y)x^sds\right|\leq F''_{N,k-1}\frac{x\log y}{(\log x)^{N+10}},\tag{4.10}\]
and 
\[\left|\frac{1}{2\pi i}\int_{\mathcal{C}(x)}\sum_{i\geq1}\phi_{i,0,N+10}(y)(s-1)^ix^sds\right|\le F'''_{N,k-1}\frac{x\log y}{(\log x)^{N+7}}.\tag {4.11}\]
\end{lemma}

\begin{proof}
For the first integral, we have \[\left|\frac{1}{2\pi i}\int_{\mathcal{C}(x)}\varepsilon'(s,y)x^sds\right|\le \sum_{0\leq i,r\leq k-1}\left|\int_{\mathcal{C}(x)}E_{i,r,k-1,N+10}(s)(\log\frac{1}{s-1})^r(\sum_{j=0}^{N'}\frac{(D_s^j\mathcal{G}_i)(1,y,-1)}{j!}(s-1)^j)x^sds\right|\]
\[\le e x (\frac{2c}{\log T}+\frac{\pi}{\log x})\sum_{0\leq i,r\leq k-1} \max_{s\in\mathcal{C}(x)}\left|E_{i,r,k-1,N+10}(s)(\log\frac{1}{s-1})^r(\sum_{j=0}^{N'}\frac{(D_s^j\mathcal{G}_i)(1,y,-1)}{j!}(s-1)^j)\right|.\]
Now, there is a constant $M_E$ such that, for all $i,r\geq0$ and $s\in\mathcal{C}(x)$, we have $|E_{i,r,k-1,N+10}(s)|\le M_E|s-1|^{N+10}$ so the above is 
\[\hspace*{-1cm}\le 4e x  k^2M_E(\frac{c}{\log T}+\frac{1}{\log x})^{N+11}(\log\log(x+1)+\ell(N)+2\pi)^{k-1}(\sum_{j=0}^{N'}\frac{|(D_s^j\mathcal{G}_i)(1,y,-1)|}{j!}(\frac{\ell(N)\log\log (x+1)+1}{\log x})^j).\tag{4.12}\]
By Theorem 3.1, there is a constant $m_{N',k-1}$, that does not depend on any parameters other than $N$, $k-1$, $Y_0$, $\epsilon$ and $\mathscr{p}$, such that for each $j\leq N'$,
       \[|(D_s^j\mathcal{G}_i)(1,y,-1)|\leq m_{N',k-1}(\log y)^{j+1}|\log\log (y+1)|^{k-1}.\]
Thus,       
\[\sum_{j=0}^{N'}\frac{|(D_s^j\mathcal{G}_i)(1,y,-1)|}{j!}(\frac{\ell(N)\log\log (x+1)+1}{\log x})^j\]
\[\leq\sum_{j=0}^{N'}\frac{m_{N',k-1}(\log y)^{j+1}(\log\log (y+1))^{k-1}}{j!}(\frac{\ell(N)\log\log (x+1)+1}{\log x})^j\]
\[\leq\log y|\log\log (y+1)|^{k-1}\sum_{j=0}^{N'}\frac{m_{N',k-1}}{j!}({\ell(N)\log\log (x+1)+1})^j\]

\[\le\log y|\log\log (y+1)|^{k-1}({\ell(N)\log\log (x+1)+1})^{N'} m_{N',k-1}(\sum_{j=0}^{N'}\frac{1}{j!}) \]

\[\le\log y|\log\log (x+1)|^{k-1}({\ell(N)\log\log (x+1)+1})^{N'} e\times m_{N',k-1}.\]
Substituting the above back into (4.12) gives that the LHS of (4.10) is 
\[\hspace*{-1.7cm}\le 4m_{N',k-1}M_Ee^2 k^2 \frac{x}{\log x^{N+11}}|\log\log(x+1)+\ell(N)+2\pi|^{k-1}({\ell(N)\log\log (x+1)+1})^{N'+N+11}\log y|\log\log (x+1)|^{k-1}.\tag{4.13}\]
Equation (4.10) directly follows from (4.13) because 
\(\frac{|\log\log(x+1)+\ell(N)+2\pi|^{k-1}({\ell(N)\log\log (x+1)+1})^{N'+N+11}|\log\log (x+1)|^{k-1}}{\log x}\)
is bounded by a constant.\\

     For the second bound, begin by exploiting the fact that the second integrand is analytic and utilize the residue theorem. This gives  \[\left|\frac{1}{2\pi i}\int_{\mathcal{C}(x)}\sum_{1\le i\le I}\phi_{i,0,N+10}(y)(s-1)^ix^sds\right|=\left|\frac{1}{2\pi i}\int_{b-\frac{i}{\log x}}^{b+\frac{i}{\log x}}\sum_{1\le i\le I}\phi_{i,0,N+10}(y)(s-1)^ix^sds\right|\]
     \[\leq \frac{x^b}{\pi\log x}\sum_{1\le i\le I}|\phi_{i,0,N+10}(y)|(\frac{c}{\log T}+\frac{1}{\log x})^i\le\frac{x({\ell(N)\log\log (x+1)+1})^{I}}{(\log x)^{\ell(N)+1}}\sum_{1\leq i\leq I}\frac{|\phi_{i,0,N+10}(y)|}{(\log x)^i}\]
     \[\le\frac{x|L_{N}\log\log (x+1)|^{I}}{(\log x)^{\ell(N)+1}}\sum_{1\leq i\leq I}\frac{|\phi_{i,0,N+10}(y)|}{(\log x)^i},\tag{4.14}\]
     for some constant $L_{N}>0$ that depends only on $N$. Now, there is a constant $L'_N$ that does not depend on any parameters other than what is stated in the lemma such that 
     \[ \frac{|L_{N}\log\log (x+1)|^{I}|\log \log (y+1)|^{k-1}}{\log x}\leq \frac{|(L_{N}+1)\log\log (x+1)|^{I+k-1}}{\log x}\le L'_N.\]
     Likewise, by (4.8), (4.9), and Theorem 3.1 there is a constant $F'''_{N,k-1}$, that does not depend on any parameters other than what is stated in the lemma, such that for each $i\leq N'$
     \[|\phi_{i,0,N+10}(y)|\le \frac{F'''_{N,k-1}}{(I+1) L'_{N}}(\log y)^{i+1}|\log \log (y+1)|^{k-1}.\tag{4.15} \]
     Equation (4.11) immediately follows by substituting (4.15) in (4.14) and noting that $\frac{\log y}{\log x} \leq 1$ as $1.9 \leq y \leq x^{\frac{1}{k}}$.
\end{proof}

\begin{remark}
Note that (4.8) and (4.9) imply that $\phi_{i,J,i}(y)=\phi_{i,J,I}(y)$ whenever $I\ge i$. In such cases, I will suppress the third index and write $\phi_{i,J}(y)=\phi_{i,J,i}(y)=\phi_{i,J,I}(y)$. Equations (4.8) and (4.9) can also be used to give an explicit formula for $\phi_{i,J,n}(y)$. In particular, \[\phi_{i,J}(y)=\sum_{l=0}^{k-1}\sum_{0\leq i'\leq i}\alpha_{i',i,J,l}(D_s^{i'}\mathcal{G}_l)(1,y,-1)\tag{4.16},\] for some constants $\{\alpha_{i',i,J,l}|0\leq i' \leq i,0\leq l \leq k-1\}.$ In the proof below, we will see that only the terms with $i\leq N$ play a relevant role and hence can have their third index suppressed. This is why the notation $\phi_{i,J}(y)$ is used in Theorem 1.1. 
\end{remark}
At this point, all that is left is to prove the main theorem.
\begin{proof}

From what we have 
\[\hspace*{-1.5cm} \frac{1}{2\pi i}\int_{\mathcal{C}(x)}\frac{(-1)^k}{k-1!}(D_z)^{k-1}f(s,y,z)\mid_{z=-1}\frac{x^s}{s}ds =\frac{1}{2\pi i}\int_{\mathcal{C}(x)}\sum_{\substack{ 1\le i\le I\\0\leq J\leq k-1} }\phi_{i,J,N+10}(y)(s-1)^i(\log{\frac{1}{s-1}})^Jx^sds  + \varepsilon_{N,k-1}'(x,y),\]
with $\varepsilon_{N,k-1}'(x,y)$ accounting for the errors considered in Lemmas 4.2 and 4.3 and any higher-order terms already controlled by the error in the statement of Theorem 1.1. By Lemma 4.2 and 4.3, $\varepsilon_{N,k-1}'(x,y)$ is controlled by the error in Theorem 1.1. Using the substitution $s=1+\frac{w}{\log x}$, this becomes 
\[=\frac{x}{2\pi i\log x}\int_{\mathcal{H}(x)}\sum_{\substack{ 1\le i\le I\\0\leq J\leq k-1} }\phi_{i,J,N+10}(y)(\frac{w}{\log x})^i(\log{\log x}-\log w)^Je^wdw  + \varepsilon_{N,k-1}'(x,y)\]
\[=\frac{x}{2\pi i\log x}\int_{\mathcal{H}(x)}\sum_{\substack{1\leq i\leq I \\1\leq J\leq k-1\\0\leq j\leq J} }\phi_{i,J,N+10}(y)(\frac{w}{\log x})^i(\log{\log x})^{J-j}(\log w)^j{J\choose j}e^wdw  + \varepsilon_{N,k-1}'(x,y),\]
where $\mathcal{H}(x)$ is the new contour in place of $\mathcal{C}(x)$ after the substitution $s=1+\frac{w}{\log x}$. At this point, note that the integrand for the term corresponding to $j=0$ is analytic in the convex hull of $\mathcal{H}(x)$. Thus,
\[\frac{x}{2\pi i\log x}\int_{\mathcal{H}(x)}\sum_{\substack{ 1\le i\le I\\0\leq J\leq k-1} }\phi_{i,J,N+10}(y)(\frac{w}{\log x})^i(\log{\log x})^{J}{J\choose 0}e^wdw \]
\[=\frac{x}{2\pi i\log x}\sum_{\substack{ 1\le i\le I\\0\leq J\leq k-1} }\phi_{i,J,N+10}(y)\frac{(\log{\log x})^{J}}{(\log x)^i}{J\choose 0}\int^{-\ell(N)\log\log (x+1)-i}_{-\ell(N)\log\log (x+1)+i}w^ie^wdw. \]
The absolute value of the above is
\[\leq\frac{x}{2\pi\log x}\sum_{\substack{1\leq i\leq I \\1\leq J\leq k-1} }|\phi_{i,J,N+10}(y)|\frac{(\log{\log x})^{J}}{(\log x)^i}{J\choose 0}\frac{2 (|\ell(N)\log\log (x+1)|+1)^i}{(\log x)^{\ell(N)}},\]
which is controlled by the error. From the above and Lemma 4.1, it follows that there is a $\varepsilon_{N,k-1}''(x,y)$, controlled by the error in the statement of Theorem 1.1, such that 
\[M_{k,\omega}(x,y)=\frac{x}{2\pi i\log x}\int_{\mathcal{H}(x)}\sum_{\substack{1\leq i\leq I \\1\leq J\leq k-1\\1\leq j\leq J} }\phi_{i,J,N+10}(y)(\frac{w}{\log x})^i(\log{\log x})^{J-j}(\log w)^j{J\choose j}e^wdw  + \varepsilon_{N,k-1}''(x,y).\]

Interchanging the sum and integral in the above and invoking Lemma 2.2 gives
\[M_{k,\omega}(x,y)=\frac{x}{\log x}\sum_{\substack{1\leq i\leq N \\1\leq J\leq k-1\\1\leq j\leq J} }\phi_{i,J}(y)\frac{\Gamma_{i,j}}{(\log x)^i}(\log{\log x})^{J-j}{J\choose j}  + \varepsilon_{N,k-1}'''(x,y),\]
where the terms with $i>N$ and the contribution from the use of Lemma 2.2 are added to $\varepsilon_{N,k-1}''(x,y)$ to result in $\varepsilon_{N,k-1}'''(x,y)$ which, in view of the error in (2.6), remains bounded by the error term in the statement of Theorem 1.1.
\medskip

Finally, to obtain the formula in Theorem 1.1, substitute the quantitative expressions for $(D_s^j\mathcal{G}_l)(1,y,-1)$ from Theorem 3.1 and get
\[=\frac{x}{\log x}\sum_{\substack{1\leq i'\leq i\leq N \\1\leq j\leq J\leq k-1\\0\leq j'\leq k-1} }\mu'_{N,k,i',i,j',j,J}(\log\log y)^{j'}(\log y)^{i'}\frac{\Gamma_{i,j}}{(\log x)^i}(\log{\log x})^{J-j}{J\choose j}+\varepsilon_{N,k-1}(x,y),\]
with (see the remark after the proof)
\[|\varepsilon_{N,k-1}(x,y)|\leq C_{N,k-1}\frac{x(\log\log (x+1))^{k-1}}{\log x}\bigg\{(\frac{\log y}{\log x})^{N+1}+\frac{1}{\log x\exp (c''\sqrt{\log y})}\bigg\}.\]
Note that in the main theorem, we denote $\mu_{N,k,i',i,j',j,J}:={J\choose j}\mu'_{N,k,i',i,j',j,J}$.

This concludes the proof of the main theorem. 
\end{proof}
\begin{remark}
The reader must be careful not to assume that all the coefficients $\mu_{N,k,i',i,j',j,J}$ in (1.2) are nonzero. In particular, note that since 
\[\frac{(-1)^k}{k-1!s}(D_z)^{k-1}f(s,y,z)\mid_{z=-1}=\sum_{\substack{ 1\le i\le I\\0\leq J\leq k-1} }\phi_{i,J,N+10}(y)(s-1)^i(\log{\frac{1}{s-1}})^J+\varepsilon_2(s,y)\]
was obtained from (4.5), all the $(\log{\frac{1}{s-1}})^J$ terms come from using $\log \zeta(s)=\log\frac{1}{s-1}+\Delta(s)$. On the other hand, all the $\log \log x$ terms come from $\log\frac{1}{s-1}$ after $s\to1+\frac{w}{\log x}$. Now notice that in (4.5) the power of $\log \zeta$ is $k-i-1$ and \[(D_s)^n\mathcal{G}_i(1,y,-1)=(D_s)^n(D_z)^ig(s,y,z)\mid_{z=-1,s=1}\ll_{n,i}(\log y)^{n+1}(\log \log (y+1))^i\] by (3.2). From here it follows that the nonzero coefficients in (1.2) satisfy $J-j+j' \le k-1$. 
\end{remark}

\begin{remark}
Note that both $\frac{1}{s\zeta(s)}$ and $\Delta(s)$ are $O(|1-s|)$ in the important contour. Hence, for $\phi_{1,J}(y)$ the derivatives $(D_s^j\mathcal{G}_i)(1,y,-1)$, with $j>0$, do not contribute and only $\Theta_{k-1-r,r,1}(s)$ and $\mathcal{G}_{k-1-r}(1,y,-1)$ do. This means,
\[\frac{(-1)^k}{k-1!}\sum_{r=0}^{k-1}{k-1\choose r}(s-1)(\log\frac{1}{s-1})^r\mathcal{G}_{k-1-r}(1,y,-1)=\sum_{0\leq J\leq k-1 }\phi_{1,J}(y)(s-1)(\log{\frac{1}{s-1}})^J,\]
and in particular
\[\phi_{1,J}(y)=\log y\sum_{0\leq j\leq k-1}\lambda_j (\log \log y)^j+O_{k-1}(\exp{(-c''\sqrt{\log y}})),\]
with $\lambda_0,\dots,\lambda_{k-1}\in\mathbb{Q}(\gamma,e^\gamma)$. Hence, $\mu_{1,k,1,1,j',j,J}\in\mathbb{Q}(\gamma,e^\gamma)$. It follows that
\[M_{k,\omega}(x,y)=\frac{x\log y}{(\log x)^2}\sum_{\substack{1\leq j\leq J\leq k-1\\0\leq j'\leq k-1} }\mu_{1,k,1,1,j',j,J}\Gamma_{1,j}(\log\log y)^{j'}(\log{\log x})^{J-j}+\varepsilon_{1,k-1}(x,y)\]
has its coefficients $\mu_{1,k,1,1,j',j,J}\Gamma_{1,j}\in\mathbb{Q}(e^\gamma,\gamma,\zeta(2),\dots,\zeta(k-2))$. I find it really interesting that the coefficients $\mu_{1,k,1,1,j',j,J}\Gamma_{1,j}$ related to the $k$-th order duality belong to $\mathbb{Q}(e^\gamma,\gamma,\zeta(2),\dots,\zeta(k-2))$, which, in view of the fact that $\zeta(2k)\in \mathbb{Q}(\pi)$, can gain transcendence degrees at most roughly half the time! 
\end{remark}

\section{Some preliminary upper bounds}
This section provides upper bounds in cases where sifting is faster. I plan to give this topic a much more thorough treatment in a subsequent paper. The familiar reader will know that in such instances, Buchstab iteration provides great insight. Indeed, this has been demonstrated in \cite{KA82} for the case $k=1$.  In fact, K. Alladi and I are currently detailing the case $k=2$ for the full range in \cite{YAKASS}. However, I will give preliminary estimates for the values of $Y_0\exp{(\mathscr{p}\frac{\log x}{\log\log (x+1)})}\leq y\leq x^{\frac{1}{k}} $. These estimates are summarized by the following results.
\begin{Prop}
Suppose that $y=Y(x)=Y_0\exp{(\mathscr{p}\frac{\log x}{\log\log (x+1)})}$ where $\mathscr{p}\in\mathbb{R}$ is some fixed power. In this case, we get the following relation
\[\tag{5.1}m_{k,\omega}(x,Y(x))\ll_k\frac{x(\log\log (x+1))^{k}}{\log x}.\]

More generally, there is a constant $\mathscr{B}_{k}$ that may depend on $k$ but no other variables or parameters, such that for all $1.9\leq y\leq x^{\frac{1}{k}}$ we have
\[\tag{5.2}{M_{k,\omega}(x,y)}\leq \mathscr{B}_{k} x \log x\log y (\log\log (x+1))^{k-1}.\]
 
 \end{Prop}

\begin{remark}
Observe that if we were to substitute $y=Y(x)=Y_0\exp{(\mathscr{p}\frac{\log x}{\log\log (x+1)})}$ in the quantitative expression of Theorem 1.1, we would get a bound that has at least one less factor of $\log\log (x+1)$.
\end{remark}

Notice that the contributions from (i) and (ii) in the adjusted contours are already controlled by Lemma 4.1  which directly applies to $y=Y(x)$. Therefore, (5.1) will follow from estimating the contribution from (iii), namely, the important contour $\mathcal{C}(x)$.

\begin{proof}[Proof of (5.1)]
As before, substituting $\log\zeta(s)=\log\frac{1}{s-1}+\Delta(s)$ in (4.1) gives 
\[\frac{(-1)^k}{k-1!}(D_z)^{k-1}f(s,y,z)\mid_{z=-1}=\frac{(-1)^k}{k-1!s\zeta(s)}\sum_{i,r\geq 0}{k-1\choose i,r,k-i-r-1}(\Delta(s))^{k-i-r-1}(\log\frac{1}{s-1})^r\mathcal{G}_i(s,y,-1).\]

I will again show that the bound (5.1) is satisfied by looking at the above finite sum term by term. For a given term, we start with 
\[\left|\frac{1}{i!r!k-i-r-1!}\int_{\mathcal{C}(x)}(\Delta(s))^{k-i-r-1}(\log\frac{1}{s-1})^r\mathcal{G}_i(s,y,-1)\frac{x^s}{s\zeta(s)}ds\right|.\]
Utilizing $|\frac{1}{s\zeta(s)}|\leq M_\zeta|1-s|$ everywhere in the important contour, the above is
\[\leq ex(\frac{2c}{\log T}+\frac{\pi}{2\log x}) M_{\zeta}(\frac{c}{\log T}+\frac{1}{\log x})(M_\Delta)^{k-i-r-1}(\log\log(x+1)+\ell(N)+2\pi)^r\max_{s\in\mathcal{C}(x)}\{|\mathcal{G}_i(s,y,-1)|\}. \tag{5.3}\]

If $y=Y(x)=Y_0\exp{(\mathscr{p}\frac{\log x}{\log\log (x+1)})}$ then, following the approach of Section 4, the bound from Theorem 3.1 gives a constant $m_{1,k-1}$ that depends only on $k-1$ such that the above is 
\[\leq ex(\frac{2c}{\log T}+\frac{\pi}{2\log x})^2m_{1,k-1}M_\zeta (M_\Delta)^{k-i-r-1}(\log\log(x+1)+\ell(N)+2\pi)^r\log{Y(x)}(\log\log (Y(x)+1))^{i}\]
\[ \ll_k\frac{x(\log\log (x+1))^{k}}{\log x},\]
where, as before, $T=\exp{(\frac{c}{\ell(N)}\frac{\log x}{\log \log (x+1)})}$ and we can choose $N=10$ so that there is no dependence on an additional variable $N$. This establishes (5.1). 
\end{proof}

All that remains is (5.2). What is interesting is that even if Theorem 3.1 is generalized to a larger range of $y$, the classical bounds on $g(s,y,-1)=-\prod_{p\leq y}\frac{1}{1-\frac{1}{p^s}}$ introduce a stubborn error term. As such, for the general sifting region, I have opted to avoid an incursion into the interior of the critical strip.

\begin{proof}[Proof of (5.2).]

First, note that, just like before, the Leibniz rule gives
\[\frac{1}{k-1!}(D_z)^{k-1}f(a,y,z)\mid_{z=1}=\frac{\zeta(a)}{k-1!}\sum_{i=0}^{k-1}{k-1\choose i}(\log\zeta(a))^{k-i-1}(D_z)^{i}g(a,y,z)\mid_{z=1}.\tag{5.4}\]
Now, as noted in the proof of Lemma 2.1, for real $\sigma>0$, we have 
\[\frac{1}{\sigma-1}\leq\zeta(\sigma)\leq \frac{\sigma}{\sigma-1},\]
and so, because $a=1+\frac{1}{\log x}$, for $x\geq 1.9$,
\[\log x\leq\zeta(a)\leq \log x +1\leq 2\log(x+1)\quad\text{and}\quad\log\log x\leq\log\zeta(a)\leq20\log\log(x+1), \]
as $\log$ here is the natural logarithm.

Substituting in equation (2.4) of Lemma 2.1 and combining this with (5.4) we get 
\[\hspace*{-1.5cm}\left|\frac{1}{2\pi i}\int_{a-iT}^{a+iT}\frac{(-1)^k}{k-1!}(D_z)^{k-1}f(s,y,z)\mid_{z=-1}\frac{x^s}{s}ds - {M_{k,\omega}(x,y)}\right|\leq 12\frac{m_{1,k-1} x\log (x+1) \log y (20\log\log (x+1))^{k-1} }{T}.
\]

In the above, I will take the non-standard choice $T=1$ and combine it with the usual $a=1+\frac{1}{\log x}$. This way, the right-hand side of the above already abides by the desired bounds. On the other hand,  
\[
   \left|\frac{1}{2\pi i}\int_{a-iT}^{a+iT}\frac{(-1)^k}{k-1!}(D_z)^{k-1}f(s,y,z)\mid_{z=-1}\frac{x^s}{s}ds\right|\leq\frac{ex}{\pi}\sum_{\substack{n=1\\p_{1}(n)> y} }^\infty \bigg|\frac{\mu(n) {\omega(n)-1\choose k-1}}{n^a}\bigg|\]
   \[=\frac{ex}{\pi}\frac{1}{k-1!}(D_z)^{k-1}f(a,y,z)\mid_{z=1}=\frac{ex}{\pi}\frac{\zeta(a)}{k-1!}\sum_{i=0}^{k-1}{k-1\choose i}(\log\zeta(a))^{k-i-1}(D_z)^{i}g(a,y,z)\mid_{z=1}\]
   \[\leq\frac{ex}{\pi}\frac{2\log(x+1)}{k-1!}\sum_{i=0}^{k-1}{k-1\choose i}(20\log\log(x+1))^{k-i-1}(D_z)^{i}g(a,y,z)\mid_{z=1}. \] 
Now (5.2) directly follows from applying the bounds of Theorem 3.1 to the expression above.   
\end{proof}
\begin{remark}
As noted at the start of this section, I plan to give the general range a thorough treatment in a future manuscript. As it stands, the bounds in this section are meant to be preliminary only and are not necessarily sharp. In fact, note that if $n\leq x$ had $p(n)>y$, then $\omega(n)<\frac{\log x}{\log y}$ and so 
\[|M_{k,\omega}(x,y)|\leq\sum_{\substack{p_{1}(n)> y\\ n\leq x} } {\omega(n)-1\choose k-1}\leq x{\frac{\log x}{\log y}-1\choose k-1}.\]
Thus, if $y>x^\frac{1}{n}$ for some integer $n>0$, then $|M_{k,\omega}(x,y)|\leq x{n-1\choose k-1}$. Likewise, if $y>\exp{(\frac{\log x}{\log\log(x+1)})}$,  then $|M_{k,\omega}(x,y)|\leq x{\log\log(x+1)-1\choose k-1}$ with $x{\log\log(x+1)-1\choose k-1}\ll_k x\log\log(x+1)^{k-1}$ which already gets rid of the $\log x \log y$ factor in (5.2).
\end{remark}
\section{A general discussion on the subradical and radical strict partial orders}
Recall that function $f$ in the variable $x$ has subradical growth if for any $\epsilon > 0$ $f(x)=o(x^{\epsilon})$ and more generally, if $f$ and $g$ are functions of $x$ then $f$ is subradically dominated by $g$ (written as $f<_{\sqrt[\forall]{\,\,}}g$) if for any $\epsilon > 0$, we have $f(x)=o(g(x)^{\epsilon})$. 

\begin{lemma}
    $f$ has subradical growth iff $f\ll\exp{g(\log{x})}$ with $g(x)=o(x)$.
\end{lemma}
\begin{proof}
 If $f\ll\exp{g(\log{x})}$ with $g(x)=o(x)$, then for any $\epsilon > 0$
 \[\lim_{x\to \infty}\frac{f(x)}{x^\epsilon}=\lim_{x\to \infty}\exp{(-\epsilon \log{x}(1-\frac{g(\log x)}{\epsilon\log x}))=0}\]
 where the last equality used $g(x)=o(x)$.
 
 Conversely, suppose that $f$ has subradical growth. Without loss of generality, we may suppose 
 that $f$ is positive (by replacing $f$ with $|f|+1$, for example). I will show that $f$ has an acceptable growth rate directly from the ``$\epsilon$-$N$" definition. Say we are given $\epsilon>0$. Since $f$ is subratical, there exists $N>1$ such that $x\geq N$ implies \[f(x)<x^\frac{\epsilon}{2}\]
 but then 
 \[0\leq \frac{\log{f(x)}}{\log x }<\frac{\epsilon}{2}\]
 whenever $x\geq N$. Since $e^x>x$, $x\geq N\implies e^x\geq N$ so 
 \[0\leq \frac{\log{f(e^x)}}{ x }<\frac{\epsilon}{2}\]
 whenever $x\geq N$. This shows that $g(x)=o(x)$ with $g(x)=\log{f(e^x)}$. The conclusion follows
\end{proof}
 Recall that for the strict partial order $<_{\sqrt[\exists]{\,\,}}$ we say that $f$ is radically dominated by $g$, and write $f<_{\sqrt[\exists]{\,\,}}g$, if there exists some $0<\epsilon < 1$ such that $f=o(g^{1-\epsilon})$. The next corollary utilizes the previous lemma to prove some claims made in the introduction.  
\begin{corollary}
    \begin{enumerate}[(i)]
    \item For any function $f$, the set $\{s:s<_{\sqrt[\forall]{\,\,}}f\}$, if nonempty, forms a ring closed under taking positive powers. That is to say, $\forall p>0$
    \[g\in\{s:s<_{\sqrt[\forall]{\,\,}}f\}\Leftrightarrow g^p\in\{s:s<_{\sqrt[\forall]{\,\,}}f\}.\]
    \item The set of subradical functions has neither a maximal element for the order $<_{\sqrt[\forall]{\,\,}}$ nor one for the order $<_{\sqrt[\exists]{\,\,}}$.
    \end{enumerate}
\end{corollary}

\begin{proof}
To prove (i), first notice that
$$s',s''\in\{s:s<_{\sqrt[\forall]{\,\,}}f\}\implies s'-s''\in\{s:s<_{\sqrt[\forall]{\,\,}}f\}.$$
Moreover, for any $p,\epsilon>0$ and $s',s''\in\{s:s<_{\sqrt[\forall]{\,\,}}f\}$, both $s'$ and $s''$ are $o(f^{\frac{\epsilon}{3(p+1)}})$ and so $$\lim_{x\to\infty}\frac{s's''}{f^\epsilon}=\lim_{x\to\infty}\frac{s'^{p}}{f^\epsilon}=0.$$
Thus, both $s's''$ and $s'^{p}$ are $o(f^{\epsilon})$ and (i) is proved.

To prove (ii), first note that from (i) it immediately follows that $<_{\sqrt[\exists]{\,\,}}$ has no maximal element as the ring of subradical functions is closed under taking positive powers. Now suppose, for the sake of contradiction, that $s$ is a maximal subradical function with respect to the order $<_{\sqrt[\forall]{\,\,}}$. By the previous lemma, $s\ll\exp{g(\log{x})}$ with $g(x)=o(x)$. Without loss of generality, we can take $g(x)$ to be positive and unbounded. Let $$s'(x)=\exp{(\sqrt{g(\log x)\log x})}=\exp{g'(\log{x})}$$ where $g'(x)=\sqrt{x g(x)}$. Now, $g'(x)=o(x)$ and hence, by the previous lemma, $s'$ is subradical. I claim that $s<_{\sqrt[\forall]{\,\,}}s'$. Indeed, for any $\epsilon>0$ we have 
$$0\leq \lim_{x\to \infty} \frac{\exp{g(\log{x})}}{(s'(x))^{\epsilon}}= \lim_{x\to \infty}\exp{(-\epsilon \sqrt{g(\log x)\log x}\bigg(1-\epsilon^{-1}\sqrt{\frac{g(\log x)}{\log x}}\bigg))}=0,$$
and so $\exp{g(\log{x})}=o(s'^{\epsilon})$. It follows that, because $s<<\exp{g(\log{x})}$, we also have $s=o(s'^{\epsilon})$ and so $s<_{\sqrt[\forall]{\,\,}}s'$. This contradicts the assumption that $s$ is maximal with respect to $<_{\sqrt[\forall]{\,\,}}$ and the proof is complete.
\end{proof}

Now, what is interesting about the above ordering is that $\frac{x}{\log x}\not<_{\sqrt[\forall]{\,\,}} x$ and, furthermore, these two functions generate the same subradical strict lower set (i.e. $\{f:f<_{\sqrt[\forall]{\,\,}} x\}=\{f:f<_{\sqrt[\forall]{\,\,}}\frac{x}{\log x}\}$). Loosely speaking, the strict partial order $<_{\sqrt[\forall]{\,\,}}$ can't distinguish between $x$ and $\frac{x}{\log x}$. However, it is also true that $<_{\sqrt[\forall]{\,\,}}$ assigns the same strict lower set to $x$ and $x^{10^{10}}$, as $\{f:f<_{\sqrt[\forall]{\,\,}} x\}=\{f:f<_{\sqrt[\forall]{\,\,}}x^{10^{10}}\}$, which certainly feels odd. It was for this reason that I considered the strict partial order $<_{\sqrt[\exists]{\,\,}}$ described in the introduction. This time, it is clear that $x$ and $x^{10^{10}}$ do not generate the same strict lower set as  $x<_{\sqrt[\exists]{\,\,}} x^{10^{10}}$. On the other hand, we still have $\{f:f<_{\sqrt[\exists]{\,\,}} x\}=\{f:f<_{\sqrt[\exists]{\,\,}}\frac{x}{\log x}\}$. That is to say, $x$ and $\frac{x}{\log x}$ still generate the same strict lower sets. The following result gives a fascinating link between the two strict partial orders and elaborates on the instance just described.
\begin{theorem}
    Two unbounded nonnegative functions $f$ and $g$ have the same $<_{\sqrt[\exists]{\,\,}}$-strict lower set iff both $\frac{f}{g}<_{\sqrt[\forall]{\,\,}} f$ and $\frac{g}{f}<_{\sqrt[\forall]{\,\,}} g$.
\end{theorem}
\begin{proof}
Suppose that ${\{s:s<_{\sqrt[\exists]{\,\,}} f\}=\{s:s<_{\sqrt[\exists]{\,\,}}g\}}$ but $\frac{f}{g}\not<_{\sqrt[\forall]{\,\,}} f$. Then, there is $\epsilon>0$ such that $f^\epsilon\ll\frac{f}{g}\Leftrightarrow  \frac{g}{f^{1-\epsilon}}\ll 1$. Moreover, the fact that $f$ is unbounded means that $\frac{g}{f^{1-\epsilon}}\ll 1=o(f^{\frac{\epsilon}{2}})$ and it follows that \[ g=o(f^{1-\frac{\epsilon}{2}})\implies f^{1-\frac{\epsilon}{2}}\not\in\{s:s<_{\sqrt[\exists]{\,\,}} g\}.\] On the other hand, clearly $f^{1-\frac{\epsilon}{2}}<_{\sqrt[\exists]{\,\,}} f \implies f^{1-\frac{\epsilon}{2}}\in\{s:s<_{\sqrt[\exists]{\,\,}} f\}$ and thus ${\{s:s<_{\sqrt[\exists]{\,\,}} f\}\not=\{s:s<_{\sqrt[\exists]{\,\,}}g\}}$ contradicting our supposition. In exactly the same way, supposing that ${\{s:s<_{\sqrt[\exists]{\,\,}} f\}=\{s:s<_{\sqrt[\exists]{\,\,}}g\}}$ and $\frac{g}{f}\not<_{\sqrt[\forall]{\,\,}} g$ leads to a contradiction. The conclusion follows.
\end{proof}

This last result ties back into the remarks made in the introduction about the significance of the sifting rate $y\leq Y_0\exp(\mathscr{p}\frac{\log x}{(\log\log (x+1))^{1+\epsilon}})$.

\begin{Prop}
The logarithm of the upper bound function $Y(x)=Y_0\exp(\mathscr{p}\frac{\log x}{(\log\log (x+1))^{1+\epsilon}})$ of the sifting is a maximal element in the set of logarithms of subradical functions. In other words, for any subradical function $s=s(x)$, either $\log s<_{\sqrt[\exists]{\,\,}}\log Y$ or $\log s$ and $\log Y$ are incomparable. 
\end{Prop}

\begin{proof}
Suppose otherwise, then $\log Y<_{\sqrt[\exists]{\,\,}}\log s\implies\frac{\log x}{(\log\log (x+1))^{1+\epsilon}}<_{\sqrt[\exists]{\,\,}}\log s$. I claim that $\log x<_{\sqrt[\exists]{\,\,}}\log s$. Indeed, there is  $0<\epsilon'<1$ such that \[\frac{\log x}{(\log\log (x+1))^{1+\epsilon}}=o((\log s)^{1-\epsilon'})\implies(\frac{\log x}{(\log\log (x+1))^{1+\epsilon}})^{^{1+\frac{\epsilon'}{2}}}=o((\log s)^{1-\frac{\epsilon'}{2}}).\] On the other hand, $\log x=o((\frac{\log x}{(\log\log (x+1))^{1+\epsilon}})^{^{1+\frac{\epsilon'}{2}}})$ and $\log x<_{\sqrt[\exists]{\,\,}}\log s$ follows. This implies that $x\ll s$, contradicting the fact that $s$ is a subradical function of $x$ and in particular $s=o(x)$.
\end{proof}

\section*{Acknowledgment}
I want to thank my doctoral advisor, Krishnaswami Alladi, for suggesting the problem, providing informed guidance, and teaching me the fundamentals of analytic number theory that made this manuscript possible. I also thank the referee for the helpful suggestions.

\end{document}